\documentclass[aop,MSNbibl,citesort,dvips]{arximspdf}
\usepackage{mathbh}
\usepackage{graphicx}

\doi{10.1214/11-AOP697}
\volume{41}
\issue{1}
\pubyear{2013}
\firstpage{329}
\lastpage{384}

\makeatletter

\renewcommand{\mid}{|}
\newcommand{\ccdott}{\cdot\, |}

\newcommand{\mmbox}[1]{\mbox{{\fontsize{8.36pt}{10pt}\selectfont{#1}}}}
\newcommand{\mmmbox}[1]{\mbox{{\fontsize{6.6pt}{8pt}\selectfont{#1}}}}

\newcommand{\ppprime}{$'$}

\newcommand{\iint}{\int\!\!\int}

\newcommand{\eqref}[1]{(\ref{#1})}

\newcommand{\dd}{\,d} 
\newcommand{\sumtwo}[2]{\mathop{\sum_{#1}}_{#2}} 
\newcommand{\er}{{\mathfrak r}}
\newcommand{\ind}{\mathbh{1}}

\newcommand{\sfrac}[2]{\frac{#1}{#2}}
\newcommand{\ssup}[1]{{{({#1})}}}
\newcommand{\tsupp}[1]{{{[{#1}]}}}

\newcommand{\cC}{\mathcal{C}}
\newcommand{\cE}{\mathcal{E}}
\newcommand{\cG}{\mathcal{G}}
\newcommand{\cI}{\mathcal{I}}
\newcommand{\cL}{\mathcal{L}}
\newcommand{\cY}{\mathcal{Y}}
\newcommand{\cZ}{\mathcal{Z}}

\newcommand{\M}{{\mathsf M}}
\newcommand{\eps}{\varepsilon}
\newcommand{\lam}{\lambda}
\newcommand{\sig}{\sigma}
\newcommand{\vphi}{\varphi}
\newcommand{\dto}{\downarrow}
\newcommand{\var}{\operatorname{var}}

\newcommand{\IP}{\mathbb{P}}
\newcommand{\IT}{\mathbb{T}}
\newcommand{\IN}{\mathbb{N}}
\newcommand{\IR}{\mathbb{R}}
\newcommand{\IE}{\mathbb{E}}
\newcommand{\iN}{\in\IN}
\newcommand{\ER}{{\mathsf{ER}}}

\newtheorem{theo}{Theorem}[section]
\newtheorem{lemma}[theo]{Lemma}
\newtheorem{prop}[theo]{Proposition}

\newproclaim{remark}[theo]{Remark}

\newproclaim{example}[theo]{Example}

\makeatother

\begin{document}
\begin{frontmatter}

\title{Random networks with sublinear preferential attachment: The giant component}
\runtitle{Preferential attachment: The giant component}

\begin{aug}
\author[A]{\fnms{Steffen} \snm{Dereich}\corref{}\ead[label=e1]{steffen.dereich@uni-muenster.de}\ead[label=u1,url]{http://www.mathematik.uni-marburg.de/\textasciitilde dereich/}}
\and
\author[B]{\fnms{Peter} \snm{M\"orters}\thanksref{t2}\ead[label=e2]{maspm@bath.ac.uk}\ead[label=u2,url]{http://people.bath.ac.uk/maspm/}}
\runauthor{S. Dereich and P. M\"orters}
\affiliation{Universit\"at M\"{u}nster and University of Bath}
\address[A]{Institut f\"{u}r Mathematische Statistik\\
Westf\"{a}lische Wilhelms-Universit\"{a}t M\"{u}nster\\
Einsteinstr. 62\\
48149 M\"{u}nster\\
Germany\\
\printead{e1}\\
\printead{u1}}
\address[B]{Department of Mathematical Sciences\\
University of Bath\\
Claverton Down\\
Bath BA2 7AY\\
United Kingdom\\
\printead{e2}\\
\printead{u2}} 
\end{aug}

\thankstext{t2}{Supported by
the EPSRC through an Advanced Research Fellowship.}

\received{\smonth{7} \syear{2010}}
\revised{\smonth{6} \syear{2011}}

%
\begin{abstract}
We study a dynamical random network model in which at every
construction step a new vertex is introduced and attached to every
existing vertex independently with a probability proportional to a
concave function $f$ of its current degree. We give a criterion for the
existence of a giant component, which is both necessary and sufficient,
and which becomes explicit when $f$ is linear. Otherwise it allows the
derivation of explicit necessary and sufficient conditions, which are
often fairly close. We give an explicit criterion to decide whether the
giant component is robust under random
removal of edges. 
We also determine asymptotically the size of the giant component and
the empirical distribution of component sizes in terms of the survival
probability and size distribution of a multitype branching random walk
associated with $f$.
\end{abstract}

%

\pdfkeywords{05C80, 60C05, 90B15, Barabasi--Albert model,
Erdos--Renyi model, power law, scale-free network,
nonlinear preferential attachment,
dynamic random graph, giant component,
cluster, multitype branching random walk}

\end{frontmatter}

\section{Introduction}

\subsection{Motivation and background}

Since the publication of the highly influential paper of Barab\'asi and
Albert~\cite{BA99} the
preferential attachment para\-digm has captured the imagination of
scientists across
the disciplines and has led to a host of, from a mathematical point of
view mostly nonrigorous,
research. The underlying idea is that the topological structure of
large networks, such as
the World-Wide-Web, social interaction or citation networks, can be explained
by the principle that these networks are built dynamically, and new
vertices prefer to be attached to vertices which have already a high
degree in the existing network.

Barab\'asi and Albert~\cite{BA99} and their followers argue that, by
building a network in which
every new vertex is attached to a number of old vertices with a
probability proportional to a linear function of
the current degree, we obtain networks whose degree distribution
follows a power law. This degree
distribution is consistent with that observed in large real networks,
but quite different from
the one encountered in the Erd\H{o}s--R\'enyi model, on which most of
the mathematical literature was
focused by this date. Soon after that, Krapivsky and Redner~\cite{KR01}
suggested to look at more
general models, in which the probability of attaching a new vertex to a
current one could be an
arbitrary function $f$ of its degree, called the \textit{attachment rule}.

In this paper we investigate the properties of preferential attachment
networks with general concave
attachment rules. There are at least two good reasons to do this: on
the one hand it turns out that global features of
the network can depend in a very subtle fashion on the function $f$,
and only the possibility to vary
this parameter gives sufficient leeway for statistical modeling and
allows a critical analysis of the robustness of
the results. On the other hand we are interested in the transitions
between different qualitative behaviors
as we pass from \textit{absence} of preferential attachment, the case of
constant attachment rules $f$, effectively corresponding
to a variant of the Erd\H{o}s--R\'enyi model, to \textit{strong forms} of
preferential attachment as given by linear attachment rules $f$.
In a previous paper~\cite{DM09} we have studied degree distributions
for such a model. We found the exact asymptotic degree
distributions, which constitute the crucial tool for comparison with
other models. The main result of~\cite{DM09}
showed the emergence of a \textit{perpetual hub}, a vertex which from
some time on remains the vertex of maximal degree, when
the tail of $f$ is sufficiently heavy to ensure convergence
of the series $\sum1/f(n)^2$. In the present paper, which is
independent of~\cite{DM09}, we look at the global connectivity features
of the network and ask for the emergence of
a \textit{giant component}, that is, a connected component comprising a
positive fraction of all vertices present.

Our first main result gives a necessary and sufficient criterion for
the existence of a giant component in terms of
the spectral radii of a family of compact linear operators associated
with $f$; see Theorem~\ref{main1}. An analysis of this result
shows that a giant component can exist for two separate reasons: \textit
{either} the tail of $f$ at infinity is sufficiently
heavy so that due to the strength of the preferential attachment
mechanism the topology of the network enforces existence
of a giant component, \textit{or} the bulk of $f$ is sufficiently large
to ensure that the edge density of the network is
high enough to connect a positive proportion of vertices. We show that
in the former case the giant component is robust under
random deletion of edges, whereas it is not in the latter case.
In Theorem~\ref{percol} we characterize the robust networks by a
completely explicit criterion.

The general approach to studying the connectivity structure in our
model is to analyze a process that systematically
explores the neighborhood of a vertex in the network. Locally this
neighborhood looks approximately like a tree,
which is constructed using a spatial branching process. The properties
of this random tree determine the connectivity
structure. We show that the asymptotic size of the giant component is
determined by the survival probability (see
Theorem~\ref{main2}), and the proportion of components with a given
size is given
by the distribution of the total number of vertices in this tree; see
Theorem~\ref{main3}.
It should be mentioned that although the tree approximation holds only
\textit{locally} it is sufficiently
powerful to give \textit{global} results through a technique called
\textit
{sprinkling}.

This approach as such is not new; for example,
it has been carried out for the class of inhomogeneous random graphs by
Bollob\'as, Janson and Riordan in
the seminal paper~\cite{BJR07}. What is new here is that the approach
is carried forward very substantially to treat
the much more complex situation of a preferential attachment model with
a wide
range of attachment functions including nonlinear ones. The increased
complexity originates in the first instance
from the fact that the presence of two potential edges in our model is
\textit{not independent} if these have the same
left end vertex. This is reflected in the fact that in the spatial
branching process underlying the construction, the offspring distributions
are not given by a Poisson process. Additionally, due to the
nonlinearity of the attachment function,
information about parent vertices has to be retained in the form of a
type chosen from an infinite-type space. Hence,
rather than being a relatively simple Galton--Watson tree, the analysis
of our neighborhoods has to be built
on an approximation by a multitype branching random walk, which
involves an infinite number of offspring and an uncountable
type space. In the light of this it is rather surprising that we are
able to get very fine explicit results, even in the fully
nonlinear case, in particular the explicit characterisation of
robustness; see Theorem~\ref{percol}. Moreover, in the nonlinear case
the abstract criterion for the existence of a giant component can be
approximated and allows explicit necessary or sufficient estimates,
which are typically rather close; see Proposition~\ref{main4}.

Although our results focus on the much harder case of nonlinear
attachment rules, they are also new in the case of linear attachment
rules $f$ and so represent very significant progress
on several fronts of research. Indeed, while the criterion for
existence of a giant component is abstract for a general attachment function,
it becomes completely explicit if this function is linear; see
Proposition~\ref{linear}. Similarly our formula for
the percolation threshold becomes explicit in the linear case, and our
result also includes behavior at criticality; see Remark~\ref{percol2}.
Fine results like this are currently unavailable
for the most studied variants of preferential attachment models with
linear attachment rules, in particular
those reviewed by Dommers et al.~\cite{DHH10}.

\subsection{The model}

We call a concave function $f \dvtx\{0,1,2,\ldots\} \longrightarrow
(0,\infty)$ with $f(0)\leq1$ and
\[
\Delta f(k):=f(k+1)-f(k)<1 \qquad\mbox{for all $k\geq0$}
\]
an \textit{attachment rule}. With any attachment rule we associate the
parameters
$\gamma^+:=\max_{k\ge0} \Delta f(k)$ and $\gamma^-:=\min_{k\ge0}
\Delta f(k)$,
which satisfy $0\leq\gamma^-\leq\gamma^+<1$. By concavity the limit
%
%
\begin{equation}\label{gammadef}
\gamma:=\lim_{n\to\infty} \frac{f(n)}{n} \mbox{ exists}\quad
\mbox{and}\quad \mbox{$\gamma =\gamma^-$.}
\end{equation}
Observe also that any attachment rule $f$
is nondecreasing with $f(k)\leq k+1$ for all $k\geq0$.


Given an attachment rule $f$, we define a growing sequence $({\mathcal
G}_N)_{N\in\IN}$
of random networks by the following iterative scheme:
\begin{itemize}
\item the network ${\mathcal G}_1$ consists of a single vertex (labeled
$1$) without edges;
\item at each time $N\geq1$, given the network ${\mathcal G}_N$, we
add a new vertex (labeled $N+1$);
\item insert for each old vertex $M$ a directed edge $N+1\to M$ with probability
\[
\frac{f(\mbox{indegree of }M\mbox{ at time }N)}{N}
\]
to obtain the network ${\mathcal G}_{N+1}$.
\end{itemize}
The new edges are inserted independently for each old vertex. Note that
our conditions on $f$
guarantee that in each evolution step the probability for adding an edge
is smaller or equal to $1$. Edges in the random network ${\mathcal
G}_N$ are \textit{dependent} if
they point toward the same vertex and independent otherwise.
Formally we are dealing with directed networks, but indeed,
by construction, all edges are pointing from the younger to the older
vertex, so that the
directions can trivially be recreated from the undirected (labeled)
graph. All the notions of
connectedness, which we discuss in this paper, are based on the \textit
{undirected} networks.

Our model differs from that studied in the majority of publications in
one respect: we do not add
a fixed number of edges in every step but a random number,
corresponding formally to the outdegree
of vertices in the directed network. It turns out (see Theorem 1.1(b)
in~\cite{DM09}) that this random number
is asymptotically Poisson distributed and therefore has very light tails.
The formal universality class of our model is therefore determined by
its asymptotic indegree distribution
which, by Theorem 1.1(a) in~\cite{DM09}, is given by the
probability weights
\[
\mu_k=\frac1{1+f(k)} \prod_{l=0}^{k-1} \frac{f(l)}{1+f(l)}
\qquad\mbox{for $k\in\IN\cup\{0\}$.}
\]
Note that these are power laws when $f(k)$ is of order $k$ (but $f$
need not be linear). More precisely,
as $k\uparrow\infty$,
\[
\frac{f(k)}{k}\to\gamma\in(0,1)\quad \Longrightarrow\quad\frac
{-{\log\mu_k}}{\log k} \to1+\frac1{\gamma},
\]
so that the LCD-model of Bollob\'as et al.~\cite{BRST01,BR03} compares
to the case $\gamma=\frac12$.

\subsection{Statement of the main results}

Fix an attachment rule $f$, and define a pure birth Markov process
$(Z_t \dvtx t \geq0)$ started in zero with generator
\[
Lg(k)= f(k) \Delta g(k),
\]
which means that the process leaves state $k$ with rate $f(k)$.
Given a suitable $0<\alpha<1$ we define a 
linear operator $A_\alpha$ on the Banach space ${\mathbf C}({\mathcal
S})$ of continuous,
bounded functions on ${\mathcal S}:=\{\ell\} \cup[0,\infty]$ with $\ell
$ being a (nonnumerical) symbol, by
\[
A_\alpha g(\tau):=\int_0^{\infty} g(t) e^{\alpha t} \,d\M(t)
+ \int_0^\infty g(\ell) e^{-\alpha t} \,d\M^\tau(t),
\]
where the increasing functions $\M$, respectively, $\M^\tau$, are given by
%
\begin{eqnarray*}
\M(t) & = & \int_0^t e^{-s} {\mathbb E}[f(Z_{s})] \,ds,\qquad
\M^\ell(t) = {\mathbb E}[Z_t], \\
\M^\tau(t) & = & {\mathbb E}[Z_t | \Delta Z_\tau=1 ] - \ind
_{[\tau,\infty)}(t) \qquad\mbox{for } \tau\in[0,\infty).
\end{eqnarray*}
We shall see in Remark~\ref{redom} that $\M^\tau\leq\M^{\tau'}$ for
all $\tau\geq\tau'\geq0$, and
therefore $\M^\infty=\lim_{\tau\to\infty} \M^\tau$ is well defined.
We shall see in Lemma~\ref{propop} that 
\[
A_\alpha1(0)<\infty\qquad\iff\qquad A_\alpha\mbox{ is a well-defined
compact operator.}
\]
In particular, the set $\mathcal I$ of parameters where $A_\alpha$ is a
well-defined (and therefore also
compact) linear operator is a (possibly empty) subinterval of $(0,1)$.

Recall that we say that \textit{a giant component exists} in the sequence
of networks $({\mathcal G}_N)_{N\in\IN}$
if the proportion of vertices in the largest connected component
${\mathcal C}_N \subset{\mathcal G}_N$ converges,
for $N\uparrow\infty$, in probability to a positive number.

\begin{theo}[(Existence of a giant component)]\label{main1}
No giant component exists if and only if there exists $0<\alpha<1$
such that $A_\alpha$ is a compact operator with spectral radius $\rho
(A_\alpha)\leq1$.
\end{theo}
\begin{example}\label{simex}
A sufficient but unnecessary criterion for existence of a giant
component is
that $\gamma\geq\frac12$, where $\gamma$ is as defined in \eqref
{gammadef}; see
Remark~\ref{operappro} below for the proof.
\end{example}

The most important example is the linear case $f(k)=\gamma k + \beta$.
In this case the family
of operators $A_\alpha$ can be analyzed explicitly; see Section \ref
{seaffinecase}. We obtain the following result.
\begin{prop}[(Existence of a giant component: linear case)]\label{linear}
If $f(k)=\gamma k + \beta$ for some $0\leq\gamma<1$ and $0<\beta\leq1$,
then there exists a giant component if and only if
\[
\gamma\geq\frac12 \quad\mbox{or}\quad \beta> \frac{(1/2-\gamma
)^2}{1-\gamma}.
\]
\end{prop}

This result corresponds to the following intuition: if the preferential
attachment is sufficiently strong (i.e., $\gamma\geq\frac12$),
then there exists a giant component in the network for purely
topological reasons and regardless of the edge density.
However, if the preferential attachment is weak (i.e.,
$\gamma<\frac12$) then a giant component exists only if the edge
density is sufficiently large.
\begin{example}
If $\gamma=0$, the model is a dynamical version of the Erd\H{o}s--R\'
enyi model sometimes called
\textit{Dubins' model}. Observe that in this case there is \textit{no}
preferential attachment. The criterion
for existence of a giant component is $\beta>\frac14$, a~fact which is
essentially known from work of
Shepp~\cite{Sh89}; see Bollob\'as, Janson and Riordan \cite
{BJR05,BJR07} for more details.
\end{example}
\begin{example}
If $\gamma=\frac12$ the model is conjectured to be in the same
universality class as the LCD-model of
Bollob\'as et al.~\cite{BRST01,BR03}. In this case we obtain that a
giant component exists regardless of the
value of $\beta$, that is, of the overall edge density. This is closely
related to the robustness of
the giant component under random removal of edges, obtained in~\cite{BR03}.
\end{example}

As the last example indicates, in some situations the giant component
is robust and survives a reduction in the edge density. To make
this precise in a general setup, we fix a parameter $0<p<1$, remove
every edge in the network independently with probability $1-p$ and
call the resulting network the \textit{percolated network}. We say the
giant component in a network is \textit{robust}, if, for every $0<p<1$,
the percolated network has a giant component.
\begin{theo}[(Percolation)]\label{percol}
Suppose $f$ is an arbitrary attachment rule and recall the definition of
the parameter $\gamma$ from \eqref{gammadef}. Then the giant component
in the
preferential attachment network with attachment rule $f$ is robust if
and only
if $\gamma\geq\frac12$.
\end{theo}
%
%
\begin{remark}\label{percol2}
The criterion $\gamma\geq\frac12$ is equivalent to the fact that the
size biased indegree distribution, with weights proportional to $k\mu
_k$, has infinite first moment.
Precise criteria for the existence of a giant component in the
percolated network
can be given in terms of the operators $(A_\alpha\dvtx\alpha\in
\mathcal I)$:
\begin{longlist}[(ii)]
\item[(i)] The giant component in the network is robust if and only if
${\mathcal I}=\varnothing$. Otherwise
the percolated network has a giant component if and only if
\[
p > \frac1{ \min_{\alpha\in\mathcal I} \rho(A_\alpha)}.
\]
\item[(ii)] In the linear case $f(k)=\gamma k + \beta$, for $\gamma>0$,
the network is robust if and only if $\gamma\geq\frac12$. Otherwise,
the percolated network has a giant component if and only if
%
%
\begin{equation}\label{percocrit}
p> \biggl( \sfrac{1}{2\gamma}-1\biggr) \Biggl( \sqrt{1+\sfrac{\gamma}{\beta}
}-1\Biggr) .
\end{equation}
\end{longlist}
%
%
Observe that running percolation with retention parameter $p$ on the
network $\cG_N$ with attachment rule $f$ leads to
a network which stochastically dominates the network with attachment
rule $pf$. Only if $f$ is constant, say $f(k)=\beta$,
these random networks coincide, and the obvious criterion for existence
of a giant component in this case is $p>\frac1{4\beta}$. This
is in line with the formal criterion obtained by letting $\gamma
\downarrow0$ in \eqref{percocrit}.
\end{remark}

We now define a multitype branching random walk, which represents an
idealization of the exploration of the neighborhood
of a vertex in the infinite network $\cG_\infty$ and which is at the
heart of our results on the sizes of connected components in the network.
A heuristic explanation of the approximation of the local neighborhoods
of typical points in the networks by this
branching random walk will be given at the beginning of Section~\ref{se7}.

In the multitype branching random walk particles have \textit{positions}
on the real line and \textit{types} in the
space ${\mathcal S}$.\setcounter{footnote}{1}\footnote
{Although the distinction of type
and space appears arbitrary at this point, it turns out that the
resulting structure of a branching random walk with a compact
typespace, rather than a
multitype branching process with noncompact typespace, is essential for
the analysis.}
The initial particle is of type $\ell$ with arbitrary starting position.
Recall the definition of the pure birth Markov process $(Z_t \dvtx t
\geq0)$.
For $\tau\geq0$, let $(Z_t^{\tsupp\tau} \dvtx t \geq0)$ be the same
process conditioned
to have a birth at time $\tau$.

Each particle of type $\ell$ in position $x$ generates offspring:
\begin{itemize}
\item to its right of type $\ell$ with relative positions at the jumps
of the process $(Z_t \dvtx t \geq0)$;
\item to its left with relative positions distributed according to the
Poisson point process $\Pi$ on $(-\infty,0]$ with
intensity measure
\[
e^{t} \IE[f(Z_{-t})] \dd t,
\]
and type being the distance to the parent particle.
\end{itemize}

%
\begin{figure}

\includegraphics{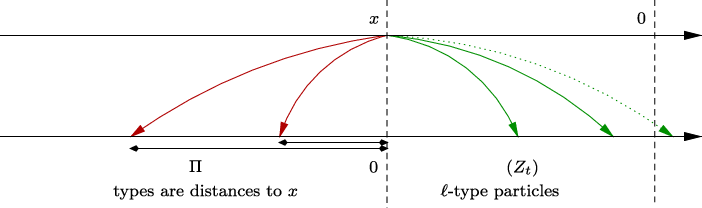}

\caption{Offspring
of an $\ell$-type particle in the branching random walk. A particle
generates finitely many offspring to its left, but infinitely many
offspring to its right.}
\label{figure1}
\end{figure}

Each particle of type $\tau\geq0$ in position $x$ generates offspring:
\begin{itemize}
\item to its left in the same manner as with a parent of type $\ell$;
\item to its right of type $\ell$ with relative positions at the jumps
of
$(Z_t^{\tsupp\tau}-\ind_{[\tau,\infty)}(t) \dvtx t \geq0)$.
\end{itemize}

%
\begin{figure}[b]

\includegraphics{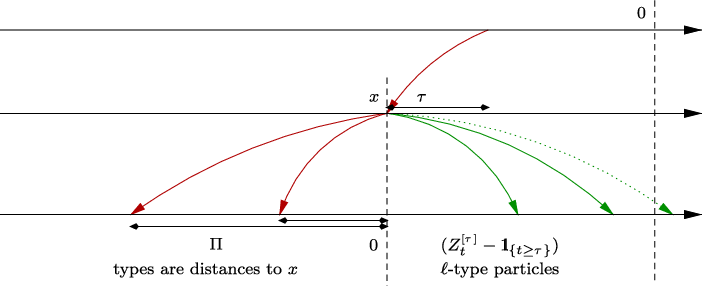}

\caption{Offspring
of a particle of type $\tau\in[0,\infty)$ in the branching random walk.
Offspring
to the right have type $\ell$, offspring to the left have type given by
the distance to the parent.}
\label{figure2}
\end{figure}

This branching random walk with infinitely many particles is called the idealized
branching random walk (IBRW); see also Figures~\ref{figure1} and~\ref{figure2}\vadjust{\goodbreak} for an illustration of the branching
mechanism.
Note
that the functions $\M$ featuring in the definition of our
operators $A_\alpha$ are derived from
the IBRW: $\M(t)$ is the expected number of particles within distance
$t$ to the left of any given particle, and
$\M^\tau(t)$ is the expected number of particles within distance $t$ to
the right of a given particle of type $\tau$.

Equally important to us is the process representing an idealization of
the exploration of the neighborhood
of a typical vertex in a large but finite network. This is the \textit
{killed} branching random walk obtained
from the IBRW by removing all particles which have a position $x> 0$
together with their entire descendancy tree.
Starting this process with one particle in position $x_0<0$ (the root),
where $-x_0$ is standard exponentially distributed,
we obtain a random rooted tree called the \textit{idealized neighborhood
tree} (INT) and denoted by $\mathfrak T$.
The genealogical structure of the tree approximates the relative
neighborhood of a typical vertex in a large but finite network. We
denote by $\# \mathfrak T$
the total number of vertices in the INT and say that the INT \textit
{survives} if this number is infinite.

The rooted tree ${\mathfrak T}$ is the weak local limit
in the sense of Benjamini and Schramm~\cite{BS01} of the sequence of
graphs in our preferential attachment model.
An interesting result about weak local limits for a different variant
of the preferential attachment network with a
linear attachment function, including the LCD-model, was recently
obtained by Berger et al.~\cite{BBCS}. In the present
paper we shall not make the abstract notion of weak local limit
explicit in our context. Instead, we go much further and give
some fine results based on our neighborhood approximation, which cannot
be obtained from weak limit theorems alone. The following
two theorems show that the INT determines the connectivity
structure of the networks in a strong sense.

\begin{theo}[(Size of the giant component)] \label{main2}
Let $f$ be an attachment rule, and denote by $p(f)$ the survival
probability of the INT. We denote by
$\cC^{\ssup1}_N$ and $\cC^{\ssup2}_N$ the largest and second largest
connected component of $\cG_N$. Then
\[
\frac{\#\cC^{\ssup1}_N}{N} \to p(f) \quad\mbox{and}\quad \frac{\#\cC^{\ssup
2}_N}{N}\to0 \qquad\mbox{in probability.}
\]
In particular, there exists a giant component if and only if $p(f)>0$.
\end{theo}

\begin{figure}

\includegraphics{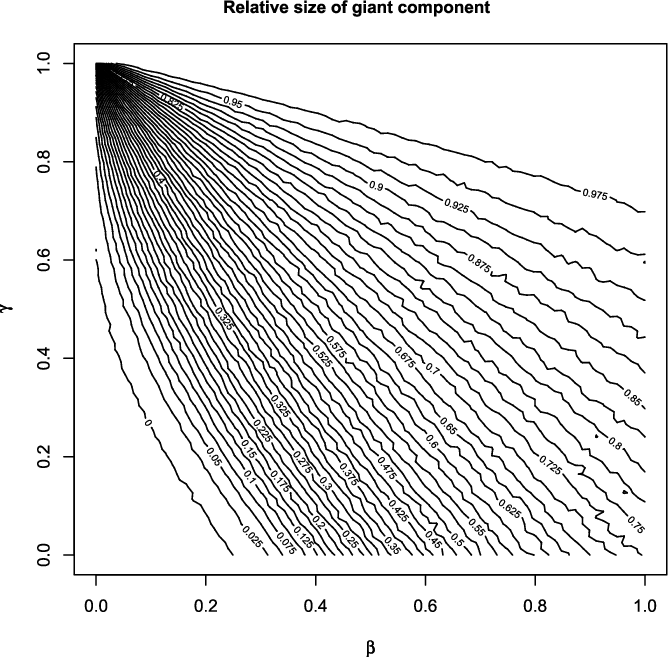}

\caption{Simulation of the proportion of vertices in the giant
component in the linear case.
The curve forming the lower envelope
is determined explicitly in Proposition \protect\ref{linear}.
The plot is based on $15\mbox{,}000$ Monte Carlo simulations of the branching
process for $80$ times $80$
gridpoints in the $(\beta,\gamma)$-plane.}
\label{figure3}\vspace*{-3pt}
\end{figure}

The results of a Monte Carlo simulation for the computation of $p(f)$ for linear~$f$ can be found
in Figure~\ref{figure3}.
The final theorem shows the cluster size distribution
in the case that no giant component exists. In this case typical
connected components, or clusters, are of finite size.
\begin{theo}[(Empirical distribution of component sizes)]\label{main3}
Let $f$ be an attachment rule,
and denote by $\cC_N(v)$ the connected component containing the
vertex $v\in\cG_N$.
Then, for every $k\in\IN$,
\[
\frac1N \sum_{v=1}^N \ind{\{\#\cC_N(v)=k\}} \longrightarrow\IP( \#
{\mathfrak T} =k )\qquad
\mbox{in probability.}
\]
\end{theo}

%

\subsection{Examples}

\subsubsection{Explicit criteria for general attachment rules}

The necessary and sufficient criterion for the existence of a
giant
component given in terms
of the spectral radius of a compact operator on an infinite-dimensional
space appears unwieldy.
%
%
\begin{figure}

\includegraphics{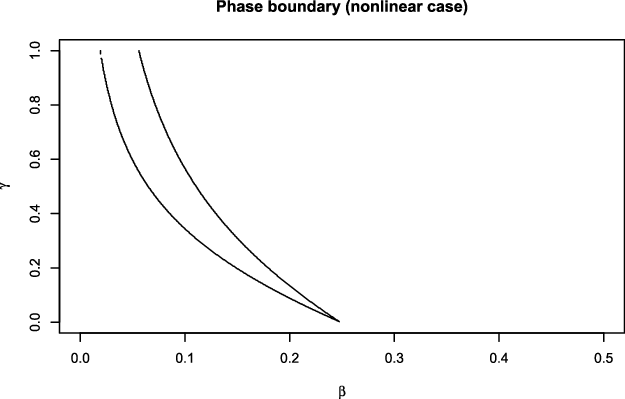}

\caption{For the attachment function $f(k)=\gamma\sqrt{k} + \beta$ the
figure shows the curves
$a[f]=\frac12$ and $a[f]+\sqrt{a[f]c[f]}=1$, which form lower and upper
bound for the boundary
between the two phases, nonexistence and existence of the giant
component, in the $(\beta,\gamma)$-plane.}
\label{figure4}\vspace*{-3pt}
\end{figure}
However, a small modification gives upper and lower bounds, which allow
very explicit necessary
or sufficient criteria that are close in many cases; see Figure
\ref{figure4}.
\begin{prop}\label{main4}
Suppose $f$ is an arbitrary attachment rule, and let
\[
a[f]:=\sum_{k=0}^\infty\prod_{j=0}^k \frac{f(j)}{1/2 + f(j)}\vadjust{\goodbreak}
\]
and
\[
c[f]:=\sum_{k=0}^\infty\prod_{j=0}^k \frac{f(j+1)}{1/2 + f(j+1)}
\geq a[f].
\]

\begin{longlist}[(ii)]
\item[(i)] If $a[f] > \frac12$, then there exists a giant
component.
\item[(ii)] If $\sfrac{1}{2} (a[f]+\sqrt{a[f]c[f]}) \leq\sfrac
{1}{2}$, then there exists no giant component.
\end{longlist}
\end{prop}
\begin{remark}\label{operappro}
\begin{itemize}
\item The term $\frac12 (a[f]+\sqrt{a[f]c[f]})$ differs from $a[f]$ by
no more than a factor~of
\[
\frac12\Biggl(1+\sqrt{\frac{1/2+f(0)}{f(0)}}\Biggr).\vadjust{\goodbreak}
\]
\item$a[f]$ converges if and only if $\gamma<\frac12$. Hence a giant
component exists if $\gamma\geq\frac12$,
as announced in Example~\ref{simex}. Otherwise there exists $\eps>0$
depending on $f(1), f(2),\ldots$
such that no giant component exists if $f(0)<\eps$.
\end{itemize}
\end{remark}
%
%
\begin{pf*}{Proof of Proposition~\ref{main4}}
(i) For a lower bound on the spectral radius we recall that $\M^\tau
\geq\M^\ell$, and therefore
we may replace $\M^\tau$ in the definition of $A_\alpha$ by $\M^\ell$.
Then $A_\alpha g(\tau)$ no
longer depends on the value of $\tau\in[0,\infty]$ but only on the fact
whether $\tau=\ell$ or otherwise.
Hence the operator collapses to become a $2 \times2$ matrix of the form
\[
\underline{A}= \pmatrix{
a(\alpha) & a(1-\alpha) \cr a(\alpha) & a(1-\alpha)}
\]
%
with
\[
a(\alpha)= \int_0^\infty e^{-\alpha t} \IE f(Z_t) \,dt.
\]
Recalling that $(Z_t \dvtx t\geq0)$ is a pure birth process with jump
rate in state $k$ given by~$f(k)$, we
can simplify this expression, using $T_k$ as the entry time into
state $k$, as follows:
\begin{eqnarray*}
\int_0^\infty e^{-\alpha t} \IE f(Z_t) \,dt
& = &\IE\sum_{k=0}^\infty f(k) \int_{T_k}^{T_{k+1}} e^{-\alpha t}
\,dt\\
& = &\sum_{k=0}^\infty f(k) \sfrac{1}{\alpha} [ \IE e^{-\alpha T_k} -
\IE e^{-\alpha T_{k+1}} ].
\end{eqnarray*}
Recalling that $T_k$ is the sum of independent exponential random
variables with parameter $f(j)$, $j=0,\ldots,k-1$,
we obtain
\[
\IE e^{-\alpha T_k} =\prod_{j=0}^{k-1} \frac{f(j)}{f(j)+\alpha}
\]
and hence
\[
a(\alpha)= \sum_{k=0}^\infty\prod_{j=0}^{k} \frac{f(j)}{f(j)+\alpha}.
\]
%
Now note that $\rho(\underline{A})=a(\alpha)+a(1-\alpha)$ and since $a$
is convex this is minimal for $\alpha=\frac12$, whence $\rho(\underline
{A}) \geq2 a(\frac12)=2 a[f]$.
This shows that the given criterion is sufficient for the existence of
a giant component.

(ii) For an upper bound on the spectral radius, we use Lemma \ref
{ledom} to see that $\M^\tau\leq\M^0$, and
therefore we may replace $\M^\tau$ in the definition of $A_\alpha$ by
$\M^0$, again reducing the operator
$A_\alpha$ to a $2 \times2$ matrix which now has the form
\[
\overline{A}= \pmatrix{
a(\alpha) & a(1-\alpha) \cr c(\alpha) & a(1-\alpha)}
\]
with $a(\alpha)$ as before and
\[
c(\alpha)= \int_0^\infty e^{-\alpha t} \IE^1 [f(Z_t)] \,dt,
\]
where $\IE^1$ is the expectation with respect to the Markov process
$(Z_t \dvtx t\geq0)$ started with
$Z_0=1$.
As before we obtain
\begin{eqnarray*}
c(\alpha)
& = &\IE^1 \biggl[\sum_{k=1}^\infty f(k) \int_{T_k}^{T_{k+1}} e^{-\alpha
t} \,dt \biggr]
= \sum_{k=1}^\infty f(k) \sfrac{1}{\alpha} \bigl[ \IE^1[ e^{-\alpha
T_k}] - \IE^1[ e^{-\alpha T_{k+1}}] \bigr]\\
& = &\sum_{k=1}^\infty f(k) \sfrac{1}{\alpha} \Biggl[ \prod_{j=2}^k \frac
{f(j-1)}{f(j-1)+\alpha}
- \prod_{j=2}^{k+1} \frac{f(j-1)}{f(j-1)+\alpha}\Biggr]\\
& = &\sum_{k=0}^\infty\prod_{j=0}^k \frac{f(j+1)}{f(j+1)+\alpha}.
\end{eqnarray*}
%
Choosing $\alpha=\frac12$, we get $\rho(\overline{A})=a[f] +\sqrt{a[f]
c[f]}$, which finishes the proof.
\end{pf*}

\subsubsection{The case of linear attachment rules}\label{seaffinecase}

We show how in the linear case $f(k)=\gamma k+\beta$ the operators
$(A_\alpha\dvtx\alpha\in\mathcal I)$ can be analyzed explicitly
and allow to infer Proposition~\ref{linear} from Theorem~\ref{main1}.
We write $\IP^k$ and $\IE^k$
for probability and expectation with respect to the Markov process
$(Z_t \dvtx t\geq0)$ started with
$Z_0=k$.
\begin{lemma}\label{linevo}
For $f(k)=\gamma k+\beta$ we have, for all $k\geq0$,
\[
\IE^k[f(Z_t)] =f(k) e^{\gamma t},\qquad
\IE^k[f(Z_t)^2] = \bigl(f(k)^2+f(k)\gamma\bigr) e^{2\gamma t}-
f(k) \gamma e^{\gamma t}
\]
and therefore
\[
d\M(t) = \beta e^{(\gamma-1)t} \,dt,\qquad d\M^{\ell
}(t) = \beta e^{\gamma t} \,dt,\qquad
d\M^\tau(t) = (\beta+\gamma)e^{\gamma t} \,dt
\]
for $\tau\in[0,\infty]$.
\end{lemma}
\begin{pf}
Recall the definition of the generator $L$ of $(Z_t \dvtx t\geq0)$.
The process
$(X_t \dvtx t\geq0)$ given by
\[
X_t=f(Z_t) - \int_0^t Lf(Z_s) \dd s=f(Z_t) - \gamma\int_0^t
f(Z_s) \dd s
\]
is a local martingale. Let $(\tau_n)_{n\iN}$ be a localizing sequence
of stopping times, and note that
\begin{eqnarray*}
\IE^k [f(Z_t)]& = &\lim_{n\to\infty} \IE^k f(Z_{t\wedge\tau_n})
=f(k)+\gamma\lim_{n\to\infty} \IE^k \int_0^{ t\wedge\tau_n} f(Z_s)
\dd s\\
&=&
f(k)+\gamma\int_0^{ t} \IE^k[ f(Z_s)] \dd s.
\end{eqnarray*}
We obtain the unique solution $\IE^k[f(Z_t)]=f(k) e^{\gamma t}$. The analogous
approach with $f$ replaced by $f^2$ gives
\begin{eqnarray*}
\IE^k [f^2(Z_t)] & = & \gamma^2 \int_0^t \IE^kf(Z_s) \,ds +
2\gamma\int_0^t \IE^k[f^2(Z_s)] \,ds +f(k)^2\\
& = & f(k) \gamma (e^{\gamma t} -1 ) + 2\gamma\int_0^{ t} \IE
^k[ f^2(Z_s)] \dd s +f(k)^2,
\end{eqnarray*}
and we obtain the unique solution
\[
\IE[f^2(Z_t)]=\bigl(f(k)^2+f(k)\gamma\bigr) e^{2\gamma t}-
f(k) \gamma e^{\gamma t}.
\]
The results for $\M$ and $\M^\ell$ follow directly from these formulas.
To characterize $\M^\tau$
for $\tau\in[0,\infty)$, we observe that, for $t\geq\tau$,
\begin{eqnarray*}
\IE[f(Z_t) \mid\Delta Z_\tau=1] &=& \sum_{k=0}^\infty\IP(Z_\tau=k)
\frac{f(k)}{\IE f(Z_\tau)} \IE^{k+1}[f(Z_{t-\tau})]\\
& = &\frac{e^{\gamma(t-2\tau)}}{\beta} \sum_{k=0}^\infty\IP(Z_\tau
=k) f(k) f(k+1)\\
& = &\frac{e^{\gamma(t-2\tau)}}{\beta} \bigl( \IE f^2(Z_{\tau}) +
\gamma \IE f(Z_{\tau}) \bigr)\\
& = &\frac{e^{\gamma(t-2\tau)}}{\beta} ( \beta^2+\beta\gamma)
e^{2\gamma\tau} = (\gamma+\beta) e^{\gamma t}
\end{eqnarray*}
and, for $t< \tau$,
\begin{eqnarray*}
\IE[f(Z_t) \mid\Delta Z_\tau=1] &=& \sum_{k=0}^\infty\IP(Z_t=k) f(k)
\frac{\IE^{k}[f(Z_{\tau-t})]}{\IE f(Z_\tau)}
\\
&=& \sum_{k=0}^\infty\IP(Z_t=k) f(k) \frac{f(k)}{f(0)}
e^{-\gamma t} = \frac{e^{-\gamma t}}{\beta} \IE[ f^2(Z_t)]\\
&=&
(\gamma+\beta) e^{\gamma t} - \gamma.
\end{eqnarray*}
From this we obtain
\[
\M^\tau(t)=\IE\bigl[ Z_t^{_{[\tau]}}\bigr] -\ind_{[\tau,\infty)}(t) =
\biggl( \sfrac{\beta}{\gamma}+1\biggr) e^{\gamma t} -1 - \sfrac{\beta}{\gamma},
\]
and, by differentiating, this implies $d\M^\tau(t)=(\beta+\gamma
)e^{\gamma t} \,dt$.
\end{pf}

\begin{pf*}{Proof of Proposition~\ref{linear}}
As $\M^\tau$ depends only on whether $\tau=\ell$ or not, the state
space ${\mathcal S}$ can be collapsed into a space
with just two points. The operator $A_\alpha$ becomes a $2 \times
2$-matrix which, as we see from the formulas below,
has finite entries if and only if
$\gamma<\alpha<1-\gamma$. This implies that there exists a giant
component if $\gamma\geq\frac12$, as in this case
the operator $A_\alpha$ is never well defined. Otherwise, denoting the
collapsed state of $[0,\infty)$ by the symbol $\er$, the matrix equals
\begin{eqnarray*}
A^{q,\er}_\alpha& = &\beta \int_0^{\infty} e^{(\gamma+\alpha-1) t}
\,dt = \frac{\beta}{1-\gamma-\alpha}\qquad \mbox{for $q\in\{\er, \ell\}
$,}\\
A^{\ell,\ell}_\alpha& = &\beta \int_0^\infty e^{(\gamma-\alpha) t}
\,dt = \frac{\beta}{\alpha-\gamma},\\
A^{\er,\ell}_\alpha& = &(\beta+\gamma) \int_0^\infty e^{(\gamma-\alpha
) t} \,dt = \frac{\beta+\gamma}{\alpha-\gamma}.
\end{eqnarray*}
Then $\rho(A_\alpha)$ is the (unique) positive solution of the
quadratic equation
\[
x^2(1-\gamma-\alpha)(\alpha-\gamma) - x(\beta-2\beta\gamma)-\beta\gamma=0.
\]
This function is minimal when the factor in front of $x^2$ is maximal,
that is, when
$\alpha=\frac12$. We note that
\[
\rho(A_{1/2})= \frac{\sqrt{\beta^2+\beta\gamma} + \beta
}{1/2 - \gamma},
\]
which indeed exceeds one if and only if
\[
\beta>\frac{(1/2-\gamma)^2}{1-\gamma}.
\]
\upqed\end{pf*}

\subsection{Overview}

The remainder of this paper is devoted to the proofs of the main results.
In Section~\ref{se2} we discuss the process describing the
indegree evolution of a fixed vertex in the network
and compare it to the process $(Z_t \dvtx t\geq0)$. The results of
this section will be frequently
referred to throughout the main parts of the proof. Section
\ref{se3} is devoted to the study of
the idealized branching random walk and explores its relation to the
properties of the family of operators \mbox{$(A_\alpha
\dvtx\alpha\in\mathcal I)$}. The main result of this section is
Lemma~\ref{opimp} which shows how survival of
the killed IBRW can be characterized in terms of these operators. Two
important tools in the proof of
Theorem~\ref{main1} are provided in Section~\ref{se4}, namely
the sprinkling argument that enables us to make
statements about the giant component from local information (see
Proposition~\ref{sprinkling}) and Lemma~\ref{initial}
which ensures by means of a soft argument that the oldest vertices are
always in large connected components.%

The core of the proof of all our theorems is provided in Sections \ref
{se1stcoupling} and~\ref{se7}. In
Section~\ref{se1stcoupling} we introduce the exploration
process, which systematically explores
the neighborhood of a given vertex in the network.
We couple this process with an analogous exploration on a random
labeled tree and show that with probability
converging to one both find the same local structure; 
see Lemma~\ref{premature}.
This random labeled tree, introduced in Section~\ref{se51}, is
still dependent on the network
size $N$, but significantly easier to study than the exploration
process itself.
Section~\ref{se7} uses further coupling arguments to relate the
random labeled tree of Section~\ref{se51} for large $N$
with the idealized branching random walk. The main result of these core
sections is summarized in Proposition~\ref{exploresurvive}.

In Section~\ref{se6} we use a coupling technique similar to
that in Section~\ref{se1stcoupling}
to produce a variance estimate for the number of vertices in components
of a given size; see Proposition~\ref{prop0825-1}. Using the machinery
provided in Sections~\ref{se4} to~\ref{se6} the proof of Theorem \ref
{main2} is completed in Section~\ref{se8} and the proof of
Theorem~\ref{main3} is completed in Section~\ref{se9}. Recall
that Theorem~\ref{main2} provides a criterion for the existence of a
giant component given in terms of the survival probability of the
killed idealized branching random walk.
In Theorem~\ref{main1} this criterion is formulated in terms of the
family of operators $(A_\alpha
\dvtx\alpha\in\mathcal I)$, and the proof of this result therefore
follows by combining Theorem~\ref{main2} with Lemma~\ref{opimp}.

The proof of the percolation result, Theorem~\ref{percol}, requires
only minor modifications of the arguments leading
to Theorem~\ref{main1} and is sketched in Section~\ref{se10}.
In a short \hyperref[app]{Appendix} we have collected some
auxiliary coupling lemmas of general nature, which are used in
Section~\ref{se7}.
Throughout the paper we use the convention that the value of positive,
finite constants $c,C$ can change from line to line,
but more important constants carry an index corresponding to the lemma
or formula line in which they were introduced.\looseness=1

\section{Properties of the degree evolution process}\label{se2}

For $m\leq n$, we denote by $\cZ[m,n]$ the indegree of vertex $m$ at
time $n$. Then, for each
$m\iN$, the degree evolution process $(\cZ[m,n] \dvtx n\geq m)$ is a
time inhomogeneous Markov process with
transition probabilities in the time-step $n\to n+1$ given by
\[
p^{\ssup n}_{k,k+1} = \frac{f(k)}n\wedge1
\quad\mbox{and}\quad p^{\ssup
n}_{k,k}=1-p^{\ssup n}_{k,k+1}\qquad
\mbox{for integers $k\geq0$.}
\]
Moreover, the evolutions $(\cZ[m, \cdot]\dvtx m\iN)$ are
independent. We suppose that under $\IP^k$
the evolution $(\cZ[m,n] \dvtx n\geq m)$ starts in $\cZ[m,m]=k$. We write
\[
P_{m,n}g(k)=\IE^k[g(\cZ[m,n])] \qquad\mbox{for any $g\dvtx\{
0,1,\ldots\}\to(0,\infty)$.}
\]
We provide several preliminary results for the process $(\cZ[m,n]
\dvtx n\geq m)$ and its continuous-time
analog $(Z_t \dvtx t\geq0)$ in this section. These form the basis for
the computations in the network.
We start by analysing the pure birth process $(Z_t \dvtx t\geq0)$
and its associated semigroup $(P_t \dvtx t \geq0)$
in Section~\ref{secont}, and then give the analogous
results for the processes $(\cZ[m,n] \dvtx n\geq m)$ in Section \ref
{sedisc}. We then compare the processes
in Section~\ref{secomp}.

\subsection{\texorpdfstring{Properties of the pure birth process $(Z_t \dvtx t\geq 0)$}
{Properties of the pure birth process (Z t: t>=0)}}\label{secont}

We start with a simple upper bound.
\begin{lemma} \label{cor0602-1}
Suppose that $f$ is an attachment rule. Then, for all $s,t\geq0$ and
integers $k\geq0$, 
\[
\IE^k[f(Z_t)]\leq f(k) e^{\gamma^+ t} \quad\mbox{and}\quad
P_{t+s}f(k)\leq e^{\gamma^+ t} P_sf(k).
\]
\end{lemma}
\begin{pf}
Note that $(Z_t \dvtx t \geq0)$ is stochastically increasing in $f$.
We can therefore obtain the result for fixed $k\geq0$
by using that $f(n)\leq f(k)+\gamma^+(n-k)$ for $n\geq k$,
and comparing with the linear model described in Lemma~\ref{linevo}.
\end{pf}

We now look at the conditioned process $(Z^{\tsupp\tau}_t \dvtx t\geq0)$.
The next two lemmas allow a comparison of the processes $(Z^{\tsupp\tau
}_t \dvtx t\geq0)$ for different values of $\tau$.
\begin{lemma}\label{comp1} 
For an attachment rule $f$, an integer $k\geq0$ and $t\geq0$, 
one has
\[
\frac{P_t f(k+1)}{P_t f(k)}\leq\frac{f(k+1)}{f(k)}
\]
for all $t\geq0$. Moreover, if $f$ is linear, then equality holds in
the display above.
\end{lemma}
\begin{pf} 
In the following, we work under the measure $\IP=\IP^{k+1}$, and we
suppose that $(U_j \dvtx j\geq0)$ is a sequence of independent
random variables, uniformly distributed in $[0,1]$, that are
independent of $(Z_t \dvtx t\geq0)$. 
We denote by $T_1,T_2,\ldots$ the random jump times of $(Z_t \dvtx
t\geq0)$ in increasing order, set $T_0=0$, and
consider the process $(Y_t \dvtx t\geq0)$ starting in $k$ that is
constant on each interval $[T_j,T_{j+1})$ and satisfies
%
%
\begin{equation}\label{eq0507-2}
Y_{T_{j+1}}=Y_{T_{j}}+\ind{\{U_j\leq f(Y_{T_j})/f(Z_{T_j})\}}.
\end{equation}
It is straightforward to verify that $(Y_t \dvtx t\geq0)$ has the
same distribution as $(Z_t \dvtx t \geq0)$ under $\IP^k$.
By the concavity of $f$ we conclude that
\[
\frac{f(Y_{T_j})}{f(Z_{T_j})}\geq\frac{f(k)+(Y_{T_j}-k) ({
f(Z_{T_j})-f(k)})/({Z_{T_j}-k})}{f(k)+(Z_{T_j}-k)
({f(Z_{T_j})-f(k)})/({Z_{T_j}-k})}
\]
and $\frac{ f(Z_{T_j})-f(k)}{Z_{T_j}-k}\leq\Delta f(k)$, so that
%
%
\begin{equation}\label{eq0507-1}
\frac{f(Y_{T_j})}{f(Z_{T_j})}\geq\frac{Y_{T_j}+{f(k)}/{\Delta
f(k)}-k}{Z_{T_j}+{f(k)}/{\Delta f(k)}-k}.
\end{equation}
Next, we couple the processes $(Y_{T_j} \dvtx j\geq0)$ and
$(Z_{T_j}\dvtx j\geq0)$ with a P\'olya urn model. Initially the urn
contains balls of two colors, blue balls of weight $B_0=\xi:=
{f(k)}/{\Delta f(k)}$, and red balls of weight one.
In each step a ball is picked with probability proportional to its
weight and a ball of the same color is inserted to the urn which
increases its weight by one.
Recalling that the total weight after $j$ draws is $j+\xi+1$, it is
straightforward to see that we can choose
the weight of the blue balls after $j$ steps as
\[
B_{j+1}=B_j +\ind{\biggl\{ U_j\leq\sfrac{B_j}{j+\xi+1}\biggr\}}.
\]
Now (\ref{eq0507-2}) and (\ref{eq0507-1}) imply that whenever we pick a
blue ball in the $j$th step, the evolution $(Y_t \dvtx t\geq0)$
increases by one at time $T_j$.
Note that $(Z_t \dvtx t \geq0)$ is independent of $(U_j \dvtx j\geq
0)$ so that
\begin{eqnarray*}\IE[Y_t| Z_t=n+k+1]-k&\geq&\IE
[B_{n}-B_0]= \frac{\xi}{1+\xi} (n+\xi+1) -\xi\\
&=&\frac{\xi n}{1+\xi}=\frac{f(k)}{f(k+1)} n
\end{eqnarray*}
and, by the concavity of $f$,
%
%
\begin{eqnarray}\label{eq0630-3}\quad
&&\IE[f(Y_t)| Z_t=n+k+1]\nonumber\\
&&\qquad\geq f(k)+ \frac{f(n+k+1)-f(k+1)}n (\IE[Y_t| Z_t=n+k+1]-k)\\
&&\qquad\geq f(k)+ \bigl(f(n+k+1)-f(k+1)\bigr) \frac{f(k)}{f(k+1)}=f(k)
\frac{f(n+k+1)}{f(k+1)},\nonumber
\end{eqnarray}
so that
\[
\frac{P_tf(k+1)}{P_tf(k)} =\frac{ \IE[ f(Z_t)]}{\IE[f(Y_t)]} \leq\frac
{f(k+1)}{f(k)}.
\]
If $f$ is linear, all inequalities above become equalities. 
\end{pf}

Next, we show that the semigroup $(P_t)$ preserves concavity.
\begin{lemma}\label{le0528-3}
For every concave and monotonically increasing $g$ and every $t\geq0$,
the function $P_t g$ is concave and monotonically increasing.
\end{lemma}
\begin{pf} 
We use an urn coupling argument similar to the one of the proof of
Lemma~\ref{comp1}. Fix $k\geq0$ and let
$(Y_t^{\ssup2} \dvtx t\geq0)$ be the pure birth process $(Z_t \dvtx
t\geq0)$ started in $Z_0=k+2$.
Denote $T_0=0$ and let $(T_j \dvtx j=1,2,\ldots)$
be the breakpoints of the process in increasing order. Suppose $(U_j
\dvtx j\geq0)$ is a sequence of independent random variables that are
uniformly distributed on $[0,1]$. For $i\in\{0,1\}$, we now denote by
$(Y_t^{\ssup i} \dvtx t\geq0)$ the step functions starting
in $k+i$ which have jumps of size one precisely at those times
$T_{j+1}$, $j\geq0$, where
\[
U_j \leq\frac{ f(Y^{\ssup i}_{T_j})}{f(Y^{\ssup2}_{T_j})}.
\]
By concavity of $f$ we get
%
%
\begin{equation}\label{eq0507-6}
\IP\bigl(\Delta Y^{\ssup1}_{T_{j+1}}=1 | \Delta Y^{\ssup0}_{T_{j+1}}=0\bigr)=
\frac{f(Y^{\ssup1}_{T_{j}}) - f(Y^{\ssup0}_{T_{j}})}
{f(Y^{\ssup2}_{T_{j}} )- f(Y^{\ssup0}_{T_{j}})}\geq\frac{Y^{\ssup
1}_{T_{j}} - Y^{\ssup0}_{T_{j}}} {Y^{\ssup2}_{T_{j}} - Y^{\ssup0}_{T_{j}}}.
\end{equation}
Let $(\bar T_j \dvtx j=1,2,\ldots)$ denote the elements of the
possibly finite set $\{ T_j\dvtx j\geq1$, $\Delta
Y^{\ssup0}_{T_{j}}=0\}$ in increasing order.
We consider a P\'olya urn model starting with one blue and one red
ball. We denote by $B_n$ the number of blue balls after $n$ steps.
By (\ref{eq0507-6}) we can couple the urn model with our indegree
evolutions such that
\[
\Delta B_j \leq\Delta Y^{\ssup1}_{\bar T_j},
\]
and such that the sequence $(B_j)_{j\iN}$ is independent of $(Y^{\ssup
2}_t \dvtx t\geq0)$ and $(Y^{\ssup0}_t\dvtx t\geq0)$.
Let $\bar g$ be the linear function on $[l,l+2+m]$ with $\bar
g(l)=g(l)$ and $\bar g(l+2+m)=g(l+2+m)$. Then
\begin{eqnarray*}
&&\IE\bigl[g\bigl(Y^{\ssup1}_t\bigr)| Y^{\ssup0}_t= l, Y^{\ssup2}_t=l+2+m
\bigr]\\
&&\qquad\geq\bar g\bigl(\IE\bigl[Y^{\ssup1}_t| Y^{\ssup0}_t= l, Y^{\ssup
2}_t=l+2+m\bigr]\bigr)\geq\bar g (l-1+\IE B_{2+m})\\
&&\qquad= \bar
g\biggl(l+1+\sfrac{m}{2}\biggr)= \sfrac{1}{2} [g(l)+g(l+2+m)].
\end{eqnarray*}
Therefore,
\[
P_tg(k+1)=\IE\bigl[g\bigl(Y^{\ssup1}_t\bigr)\bigr]\geq\tfrac{1}{2} \bigl[\IE\bigl[g\bigl(Y^{\ssup
0}_t\bigr)\bigr]+ \IE\bigl[g\bigl(Y^{\ssup2}_t\bigr)\bigr]\bigr]
=\tfrac{1}{2} [ P_tg(k)+P_tg(k+2)],
\]
which implies the concavity of $P_tg$.
\end{pf}

The fact that the semigroup preserves concavity allows us to generalize
Lem\-ma~\ref{comp1}.
\begin{lemma}\label{comp3}
For an attachment rule $f$ and integers $k\geq0$ and $s,t\geq0$, one has
\[
\frac{P_{t+s}f(k+1)}{P_{t+s} f(k)} \leq\frac{P_{s}f(k+1)}{P_{s} f(k)}.
\]
\end{lemma}
\begin{pf}
The statement follows by a slight modification of Lemma~\ref{comp1}. We
use $Z$ and $Y$ as in the proof of the latter lemma and observe that by
Lemma~\ref{le0528-3} the function
\[
g(k):=P_s f(k)
\]
is concave and increasing. Similarly as in (\ref{eq0630-3}) we get
\begin{eqnarray*}
&&\IE[g(Y_t)|Z_t=n+k+1]\\
&&\qquad\geq g(k)+ \frac{g(n+k+1)-g(k+1)}{n} (\IE[Y_t|Z_t=n+k+1]-k)\\
&&\qquad\geq g(k)+ \bigl(g(n+k+1)-g(k+1)\bigr) \frac{f(k)}{f(k+1)}\\
&&\qquad\geq g(k)+ \bigl(g(n+k+1)-g(k+1)\bigr) \frac{g(k)}{g(k+1)}=g(n+k+1)\frac
{g(k)}{g(k+1)}.
\end{eqnarray*}
The rest of the proof is in line with the proof of Lemma~\ref{comp1}.
\end{pf}
\begin{lemma}[(Stochastic domination)] \label{ledom} One can couple the
process $(Z^{\tsupp\tau}_t \dvtx\break t\geq0)$ with start in
$Z^{\tsupp\tau}_0=k$
and the process $(Z_t \dvtx t\geq0)$ with start in $Z_0=k+1$ in such
a way that
\[
\bigl\{t>0 \dvtx\Delta Z^{\tsupp\tau}_t=1\bigr\}\subset\{t>0 \dvtx\Delta
Z_t=1\}\cup\{\tau\}.
\]
In particular, this implies that $Z^{\tsupp\tau}_t+\ind{\{t<\tau\}}\leq
Z_t$ for all $t\geq0$.
In the linear case we have equality in both formulas.
\end{lemma}
\begin{pf} 
Suppose $(Y^{\ssup2}_t \dvtx t\geq0)$ has the distribution of $(Z_t
\dvtx t\geq0)$ with start in $Z_0=k+1$, let $T_0=0$
and $(T_j \dvtx j=1,2,\ldots)$ the times of discontinuities of
$(Y^{\ssup2}_t \dvtx t\geq0)$ in increasing order. Denote
by $(U_j \dvtx j\geq0)$ a sequence of independent random variables\vadjust{\goodbreak}
that are uniformly distributed on $[0,1]$. Now define
$(Y^{\ssup1}_t \dvtx t\geq0)$ as the step function starting in $k$
which increases by one: (i) at time $T_{j+1}<\tau$ if
%
%
\begin{equation}\label{eq0507-4}
U_j\leq\frac{f(Y^{\ssup1}_{T_{j}})}{f(Y^{\ssup2}_{T_{j}})}
\frac{P_{\tau-T_{j+1}}f(Y^{\ssup
1}_{T_{j}}+1)}{P_{\tau-T_{j+1}}f(Y^{\ssup1}_{T_{j}})},
\end{equation}
(ii) at time $\tau$ and (iii) at time $T_{j+1}>\tau$ if
%
%
\begin{equation}\label{eq0507-5}
U_j\leq\frac{f(Y^{\ssup1}_{T_{j}\vee\tau})}{f(Y^{\ssup2}_{T_{j}})} .
\end{equation}
Clearly, we have $Y^{\ssup1}_t+1\leq Y^{\ssup2}_t$ for all $t\in
[0,\tau)$ and $Y^{\ssup1}_t\leq Y^{\ssup2}_t$ for general $t\geq0$.
Moreover, by Lemma~\ref{comp1}, the right-hand sides of
inequalities (\ref{eq0507-4}) and (\ref{eq0507-5}) are not greater than
one, and it is straightforward to verify that
$(Y^{\ssup1}_t \dvtx t\geq0)$
has the same law as the process $(Z^{\tsupp\tau}_t \dvtx t\geq0)$
with start in $Z^{\tsupp\tau}_0=k$.
\end{pf}
\begin{remark}\label{redom} In analogy to above, one can use Lemma \ref
{comp3} to prove that two evolutions $Z^\tsupp{\sig}$ and
$Z^\tsupp{\tau}$ started in $k$ with $0< \sig\leq\tau$ can be coupled
such that
\[
\bigl\{t\geq0\dvtx \Delta Z^{[\tau]}=1\bigr\}\setminus\{\tau\}
\subset\bigl\{t\geq0\dvtx
\Delta Z^{[\sigma]}=1
\bigr\}\setminus\{\sig\}.
\]
%
\end{remark}

\subsection{\texorpdfstring{Properties of the degree evolutions $(\cZ[m,n] \dvtx n\geq m)$}
{Properties of the degree evolutions (Z[m,n]: n>=m)}}
\label{sedisc}

For the processes $(\cZ[m,n] \dvtx n\geq m)$ we get an analogous
version of Lemma~\ref{cor0602-1}.
\begin{lemma}\label{le0602-1}
For any attachment rule $f$, and all integers $k\geq0$ and $0<m\leq n$,
\[
\IE^k[f(\cZ[m,n])] \leq f(k) \biggl(\frac{n}m\biggr)^{\gamma^+}.
\]
\end{lemma}
\begin{pf}
Note that $(Y_n \dvtx n\geq m)$
with $Y_n:=f(\cZ[m,n]) \prod_{i=m}^{n-1} (1+\frac{\gamma^+}i)^{-1}$ is
a supermartingale. Hence
\[
\IE^k [f(\cZ[m,n])] \leq f(k) \prod_{i=m}^{n-1} \biggl(1+\sfrac
{\gamma^+}{i}\biggr)
\leq f(k) \biggl(\frac{n}m\biggr)^{\gamma^+}.
\]
\upqed\end{pf}

We also get the following analog of Lemma~\ref{comp1}.
\begin{lemma}\label{comp2}
For an attachment rule $f$ and integers $k\geq0$ and $0<m\leq n$ one has
\[
\frac{P_{m,n}f(k+1)}{P_{m,n}f(k)}\leq\frac{f(k+1)}{f(k)}.
\]
If $f$ is linear and $f(k+1+l)\leq m+l$ for all $l\in\{0,\ldots,n-m-1\}
$, then equality holds.\vadjust{\goodbreak}
\end{lemma}
\begin{pf}
The statement follows by a slight modification of the proof of
Lem\-ma~\ref{comp1}.\vspace*{-2pt}
\end{pf}

We now provide two lemmas on stochastic domination of the degree evolutions.\vspace*{-2pt}
%
\begin{lemma}[(Stochastic domination I)]\label{extra}
For any integers $0<m\le n_1<\cdots<n_j$ the process $(\cZ[m,n] \dvtx
n\geq m)$ conditioned on the event
$\Delta\cZ[m,n_i]=0$ for all $i\in\{1,\ldots,j\}$ is stochastically
dominated by the unconditional process.\vspace*{-2pt}
\end{lemma}
\begin{pf}
First suppose that $m<n_1$. For any $k\geq0$, we have
\begin{eqnarray*}
&&\IP^k( \Delta\cZ[m,m]=1 | \Delta\cZ[m,n_i]=0\ \forall i\in\{
1,\ldots,j\} )\\
&&\qquad= \frac{f(k)}{m} \frac{\IP^{k+1}(\Delta\cZ[m+1,n_i]=0\ \forall
i)}{\IP^{k}(\Delta\cZ[m,n_i]=0\ \forall i)}.
\end{eqnarray*}
The denominator on the right is equal to
\begin{eqnarray*}
&&\sfrac{f(k)}{m} \IP^{k+1}(\Delta\cZ[m+1,n_i]=0\ \forall i) + \biggl(1-
\sfrac{f(k)}{m}\biggr) \IP^{k}(\Delta\cZ[m+1,n_i]=0\ \forall i)\\
&&\qquad\geq\IP^{k+1}(\Delta\cZ[m+1,n_i]=0\ \forall i),
\end{eqnarray*}
and hence we get
%
%
\begin{eqnarray}\label{quey}
&&\IP^k( \Delta\cZ[m,m]=1 | \Delta\cZ[m,n_i]=0\ \forall i\in\{
1,\ldots,j\} )\nonumber\\[-8pt]\\[-8pt]
&&\qquad\leq\frac{f(k)}{m} = \IP^k( \Delta\cZ[m,m]=1 ),\nonumber
\end{eqnarray}
which is certainly also true if $m=n_1$. The result follows by induction.\vspace*{-2pt}
\end{pf}

The next lemma is the analog of 
Lemma~\ref{ledom}.\vspace*{-2pt}
%
%
\begin{lemma}[(Stochastic domination II)]\label{ledom2}
For integers $0 \leq k < m < n$ there exists a coupling of the process
$(\cZ[m,l] \dvtx l\geq m)$ started in $\cZ[m,m]=k$
and conditioned on $\Delta\cZ[m,n]=1$ and the unconditional process
$(\cZ[m,l] \dvtx l\geq m)$ started in $\cZ[m,m]=k+1$
such that for the coupled random evolutions, say $(\cY^{\ssup1}[l]
\dvtx l\geq m)$ and $(\cY^{\ssup2}[l] \dvtx l\geq m)$, one has
\[
\Delta\cY^{\ssup1}[l]\leq\Delta\cY^{\ssup2}[l] +\ind{\{l=n\}},
\]
and therefore in particular $\cY^{\ssup1}[l]\leq\cY^{\ssup2}[l]$ for
all $l\geq m$.\vspace*{-2pt}
\end{lemma}
\begin{pf}
Note that
\begin{eqnarray*}
\IP^k (\Delta \cZ[m,m]=1| \Delta\cZ[m,n]=1)&=& \frac{\IP^k (\Delta
\cZ[m,m]=1, \Delta\cZ[m,n]=1)}{\IP^k (\Delta\cZ[m,n]=1)}\\
&=& \frac{({f(k)}/m) \IE^{k+1}[f(\cZ[m+1,n])](1/n)}{\IE^k[ f (\cZ
[m,n])](1/n)}\\
&=&\frac{f(k)}m \frac{P_{m+1,n}f(k+1)}{P_{m,n} f(k)}.
\end{eqnarray*}
By Lemma~\ref{comp2}, we get
\[
\IP^k (\Delta\cZ[m,m]=1| \Delta\cZ[m,n]=1)\leq\frac{f(k)}m
\frac{P_{m+1,n}f(k+1)}{P_{m+1,n} f(k)} \leq\frac{f(k+1)}{m}.
\]
Now the coupling of the processes can be established as in Lemma~\ref{ledom}.
\end{pf}
\begin{lemma}\label{le0416-1}
For all $m\leq n\leq n'$ one has
\[
\IP(\Delta\cZ[m,n]=1)\geq\IP(\Delta\cZ[m,n']=1).
\]
\end{lemma}
\begin{pf}It suffices to prove the statement for $n'=n+1$ and $n\geq m$
arbitrary.
The statement follows immediately from
\[
\IP(\Delta\cZ[m,n]=1)= \frac1n \IE[ f(\cZ[m,n])]=\frac1n \sum
_{k=0}^\infty\IP(\cZ[m,n]=k) f(k)
\]
and
\begin{eqnarray*}
&&\IP(\Delta\cZ[m,n+1]=1) \\ 
&&\qquad= \frac{1}{n+1} \sum_{k=0}^\infty\IP(\cZ[m,n]=k) \biggl[ \frac
{f(k)}{n} f(k+1) + \biggl(1-\frac{f(k)}{n}\biggr) f(k)\biggr]\\
&&\qquad=\frac1n \sum_{k=0}^\infty\frac{n+\Delta f(k) }{n+1} f(k) \IP
(\cZ
[m,n]=k).
\end{eqnarray*}
\upqed\end{pf}

We finally look at degree evolutions $(\cZ[m,n] \dvtx n\geq m)$
conditioned on both the existence and nonexistence of
some edges. In this case we cannot prove stochastic domination, and
comparison requires a constant factor.
\begin{lemma}\label{condaway}
Suppose that $(c_N)_{N\iN}, (n_N)_{N\iN}$ are sequences of integers
such that $\lim_{N\to\infty} n_N=\infty$ and
$c_N^2 n_N^{\gamma^+-1}$ is bounded from above. Then there exists a
constant $C_{\mmbox{\ref{condaway}}}>0$,
such that for all ${\mathcal I}_0, {\mathcal I}_1$ disjoint subsets of
$\{{n_N}, \ldots, N\}$ with $\#{\mathcal I}_0\leq c_N$ and $\#{\mathcal
I}_1\leq1$ and, for any $m\in\{1,\ldots,N\}$ with $n\geq m$, we have
\begin{eqnarray*}
&&\IP( \Delta\cZ[m,n-1]=1 | \Delta\cZ[m,i]=1\ \forall i\in\cI_1,
\Delta\cZ[m,i]=0\ \forall i\in\cI_0)\\
&&\qquad\leq C_{\mmbox{\ref{condaway}}} \IP( \Delta\cZ[m,n-1]=1 |
\Delta
\cZ[m,i]=1\ \forall i\in\cI_1).
\end{eqnarray*}
\end{lemma}
\begin{pf}
We have
\begin{eqnarray*}
&&\IP( \Delta\cZ[m,n-1]=1 | \Delta\cZ[m,i]=1\ \forall i\in\cI_1,
\Delta\cZ[m,i]=0\ \forall i\in\cI_0)\\
&&\qquad\leq\frac{\IP( \Delta\cZ[m,n-1]=1 | \Delta\cZ[m,i]=1\ \forall
i\in
\cI_1)}
{\IP(\Delta\cZ[m,i]=0\ \forall i\in\cI_0 | \Delta\cZ[m,i]=1\ \forall
i\in\cI_1)},
\end{eqnarray*}
and it remains to bound the denominator from below by a positive constant.

Using Lemma~\ref{ledom2} and denoting $k=\#\cI_1$ we obtain that
\begin{eqnarray*}
&&\IP(\Delta\cZ[m,i]=0\ \forall i\in\cI_0 | \Delta\cZ
[m,i]=1\ \forall i\in\cI_1) \\
&&\qquad\geq\IP^1(\Delta\cZ[m,i]=0\ \forall i\in\cI_0) \geq\prod
_{j\in\cI_0} \IP^1(\Delta\cZ[m,j]=0)\\
&&\qquad= \prod_{j\in\cI_0} \biggl\{1-\frac{\IE^1[f(\cZ[m,j])]}{j}\biggr\}.
\end{eqnarray*}
By Lemma~\ref{le0602-1} the expectation is bounded from above by
$f(k)j^{\gamma^+}$ and
moreover $f(k)\leq k+1 \leq2 c_N$ for $N$ large enough. Hence we get,
\[
\prod_{j\in\cI_0} \biggl\{1-\frac{\IE^1[f(\cZ[m,j])]}{j}\biggr\}
\geq\prod_{j\in\cI_0} \{1- 2 c_N j^{\gamma^+-1}\}
\geq(1- 2 c_N {n_N}^{\gamma^+-1})^{c_N}
\]
using that $\# \cI_0\leq c_N$. As $c_N^2 {n_N}^{\gamma^+-1}$ is bounded from
above, the expression on the right is bounded from zero.
This implies the statement.
\end{pf}

\subsection{Comparing the degree evolution and the pure birth process}
\label{secomp}

The aim of this section is to show that the processes $(\cZ[m,n] \dvtx
n\geq m)$ and
$(Z_t \dvtx t \geq0)$ are intimately related. To this end, we set
%
%
\begin{equation}\label{tndef}
t_n:=\sum_{k=1}^{n-1} \frac1k \quad\mbox{and}\quad \Delta
t_n:=t_{n+1}-t_n=\frac1n.
\end{equation}

\begin{lemma}\label{le0525-1} 
For fixed $n\in\IN$, one can couple the random variables $Z_{\Delta t_{n}}$
and $\cZ[n,n+1]$ under $\IP^k$
such that, almost surely,
%
\[
\IP(Z_{\Delta t_{n}}\not=\cZ[n,n+1])\leq\bigl( f(k+1) \Delta t_n\bigr)^2
\quad\mbox{and}\quad (k+1)\wedge Z_{\Delta t_{n}} \leq\cZ[n,n+1].
\]
\end{lemma}
\begin{pf}
Note that
\begin{eqnarray*}
\IP^k(Z_{\Delta t_{n}}=k+1)&=& f(k) \Delta t_n e^{-f(k) \Delta t_n}
\frac
1{\Delta t_n} \int_0^{\Delta t_n} e^{-\Delta f(k) u} \dd u\\
&\geq& f(k) \Delta t_n e^{-f(k+1) \Delta t_n}.
\end{eqnarray*}
The same lower bound is valid for the probability $\IP^{k}(\cZ
[n,n+1]=k+1)$. Moreover,
\[
\IP^k(Z_{\Delta t_{n}}=k)= e^{-f(k) \Delta t_n}\geq\bigl(1-f(k) \Delta
t_n\bigr)\vee0= \IP^{k}( \cZ[n,n+1]=k).\vadjust{\goodbreak}
\]
Hence, we can couple $Z_{\Delta t_{n}}$ and $\cZ[n,n+1]$ under $\IP
^{k}$ such that that they differ with probability less than
%
%
\begin{eqnarray}\label{eqcoupling}
&&1-\bigl[f(k) \Delta t_n e^{-f(k+1) \Delta t_n}+1-f(k) \Delta
t_n\bigr]\nonumber\\[-8pt]\\[-8pt]
&&\qquad=f(k) \Delta t_n \bigl(1- e^{-f(k+1)\Delta t_n }\bigr)
\leq\bigl( f(k+1) \Delta
t_n\bigr)^2,\nonumber
\end{eqnarray}
and moreover we have $ (k+1)\wedge Z_{\Delta t_{n}}\leq\cZ[n,n+1]$.
\end{pf}
\begin{prop}\label{SGapprox} There exist constants $n_0\iN$ and
$C_{\mmbox{\ref{SGapprox}}}>0$
such that for all integers $n_0\leq m\leq n$ and $0\leq k<m$,
\[
|P_{m,n} f(k)-P_{t_n-t_m} f(k)| \leq C_{\mmbox{\ref{SGapprox}}} \frac
{f(k)}{m} P_{m,n} f(k).
\]
\end{prop}

The proof of the proposition uses several preliminary results on the
semigroups $(P_t \dvtx t \geq0)$ and
$(P_{m,n} \dvtx n\geq m)$, which we derive first.
For a stochastic domination argument we introduce a further time
inhomogeneous Markov process. For
integers $n,k\geq0$, we suppose that
\begin{eqnarray*}
\tilde\IP^k (\cZ[n,n+1]=k+1) &=& 1-\tilde\IP^k (\cZ[n,n+1]=k)\\
&=& \biggl(\frac{f(k)}n+\frac12 f(k) \Delta f(0) e^{\Delta f(0)} \frac
1{n^2}\biggr)\wedge1.
\end{eqnarray*}
The corresponding semigroup is denoted by $(\tilde P_{m,n})_{m\leq n}$.
\begin{lemma}\label{le0528-2}
Assume that there exists $n_0\iN$ such that, for all integers \mbox{$n\geq n_0$},
%
%
\begin{equation}\label{eq0528-4}
\frac{f(n)}n+\frac12 f(n) \Delta f(0) e^{\Delta f(0)} \frac
1{n^2}\leq1.
\end{equation}
Then, for all integers $n\geq n_0$ and $0\leq k\leq n$, and an
increasing concave $g\dvtx\{0,1,2,\ldots\} \to\IR$,
\[
P_{\Delta t_n}g(k) \leq\tilde P_{n,n+1}g(k) .
\]
\end{lemma}
\begin{pf} Consider $\bar f(l)= f(k) +\Delta f(k) (l-k)$.
Note that by 
comparison with the linear model
\[
f(k)+\Delta f(k) (\IE^k[ Z_t]-k) = \IE^k[\bar f(Z_t)]\leq f(k)
e^{\Delta f(k) t}.
\]
Hence, for $t\in[0,1]$, using that $e^x\leq1+x+ \frac12 x^2 e^x$ for
$x\geq0$,
\[
\IE^k[ Z_t]-k\leq\frac{f(k)}{\Delta f(k)} \bigl(e^{\Delta f(k) t}-1\bigr)\leq
f(k) t+\sfrac{1}{2} f(k) \Delta f(k) e^{\Delta f(k) t} t^2.
\]
Therefore, $\IE^k[ Z_{\Delta t_n}]\leq\tilde\IE^k[\cZ[n,n+1]]$ for
all $n\geq n_0$.
As $g$ is increasing and concave, and $\cZ$ has only increments of size
one, we get
\begin{eqnarray*}
\IE^k[ g(Z_{\Delta t_n})] &\leq& g(k) + \bigl( g(k+1)-g(k)
\bigr) \IE^k[ Z_{\Delta t_n}-k]\\
&\leq& g(k) + \bigl( g(k+1)-g(k) \bigr) \tilde\IE^k\bigl[ \cZ[n,n+1]-k \bigr]
\\
&=&\tilde\IE^k\bigl[g(\cZ[n,n+1])\bigr]
\end{eqnarray*}
as required to complete the proof.
\end{pf}
\begin{lemma}\label{le0528-1}
There exists a constant $C_{\mmbox{\ref{le0528-1}}}>0$, depending on
$f$, such
that for all integers $0\leq k \leq m$ and
$0<m\leq n$, we have
\[
\tilde P_{m,n} f(k) \leq C_{\mmbox{\ref{le0528-1}}} P_{m,n} f (k).
\]
\end{lemma}
\begin{pf} For $n,m\iN$ with $n\geq m$ let $c_{m,n}:=\prod_{l=m}^{n-1}
(1+\frac\kappa{l^2})$ where
$\kappa:=\frac12 (\Delta f(0))^2 e^{\Delta f(0)}$. We prove by
induction (over $n-m$) that for all $0<m\leq n$ and
$0\leq k\leq m$,
\[
\tilde P_{m,n} f(k) \leq c_{m,n} P_{m,n} f (k).
\]
Certainly the statement is true if $n=m$. Moreover, we have
\[
\tilde P_{m,n+1} f(k)= P_{m,{m+1}} \tilde P_{m+1,n+1} f(k) + (\tilde
P_{m,m+1}-P_{m,m+1}) \tilde P_{m+1,n+1} f(k),
\]
and applying the induction hypothesis we get
\[
\tilde P_{m,n+1} f(k)\leq c_{m+1,n+1} P_{m,{n+1}} f(k) + (\tilde
P_{m,m+1}-P_{m,m+1}) \tilde P_{m+1,n+1} f(k).
\]
Moreover, for a function $g \dvtx\{0,1,2,\ldots\}\to\IR$, we have
%
%
\begin{equation}\label{eq0528-5}
(\tilde P_{m,m+1}-P_{m,m+1}) g(k)\leq\frac12 f(k) \Delta f(0)
e^{\Delta f(0)} \frac1{m^2} \Delta g(k).
\end{equation}
Note that the transition probabilities of the new inhomogeneous Markov
process have a particular product structure:
for all integers $a\geq1$ and $b\geq0$, one has
\[
\tilde\IP^b(\cZ[a,a+1]=b+1)=\bigl(\psi_a\cdot f(b)\bigr)\wedge1\qquad
\mbox{for }\psi_a:=\frac{1}a+\frac12 \Delta f(0) e^{\Delta f(0)}
\frac1{a^2}.
\]
This structure allows one to literally translate the proof of Lemma \ref
{comp2} and to obtain
\[
\frac{\tilde P_{a_1,a_2}f(b_2)}{\tilde P_{a_1,a_2}f(b_1)} \leq\frac
{f(b_2)}{f(b_1)}
\]
for integers $a_1, a_2\geq1$ and $b_1, b_2\geq0$ with $a_1\leq a_2$
and $b_1\leq b_2$.
Consequently, using (\ref{eq0528-5}) and the induction hypothesis,
%
%
\begin{eqnarray}\label{eq0423-3}
&&(\tilde P_{m,m+1}-P_{m,m+1}) \tilde P_{m+1,n+1} f(k)\nonumber\\
&&\qquad\leq\frac12 f(k) \Delta f(0) e^{\Delta f(0) } \frac1{m^2}
\frac{\Delta f(k)}{f(k)}\tilde P_{m+1,n+1} f(k)\\
&&\qquad\leq\frac{\kappa}{m^2} \tilde P_{m+1,n+1} f(k)\leq\frac\kappa{m^2}
c_{m+1,n+1} P_{m+1,n+1}
f(k).\nonumber
\end{eqnarray}
Altogether, we get
\[
\tilde P_{m,n+1} f(k)\leq\biggl(1+\frac{\kappa}{m^2}\biggr) c_{m+1,n+1}
P_{m,{n+1}} f(k)=c_{m,n+1}P_{m,{n+1}} f(k),
\]
and the statement follows since all constants are uniformly bounded by
$\prod_{l=1}^\infty(1+\frac{\kappa}{l^2})<\infty$.
\end{pf}
\begin{pf*}{Proof of Proposition~\ref{SGapprox}} 
We choose $n_0$ as in Lemma~\ref{le0528-2}, and let $k,m,n$ be integers
with $n_0\leq m\leq n$ and $0\leq k\leq m$.
We represent $\IE^k[f(\cZ[m$, $n])]-\IE^k[f(Z_{t_n-t_m})]$ as the
telescoping sum
%
%
\begin{equation}\label{eq0423-2}\quad
P_{m,n} f(k)- P_{t_n-t_m} f(k)= \sum_{l=m}^{n-1} \underbrace{P_{m,l}
(P_{l,l+1}-P_{t_{l+1}-t_l}) P_{t_n-t_{l+1}} f(k)}_{=:\Sigma_l}.
\end{equation}
In the following, we fix $l\in\{m,\ldots,n-1\}$ and analyze the summand
$\Sigma_l$.
First note that by Lemma~\ref{comp1}, one has for arbitrary integers
$0\leq a\leq b$,
%
%
\begin{eqnarray}\label{eq0506-1}
\vphi(a,b) :\!&= & \IE^{b} [f (Z_{t_n-t_{l+1}})]-\IE^{a} [f
(Z_{t_n-t_{l+1}})]\nonumber\\[-8pt]\\[-8pt]
&\leq& \frac{f(b)-f(a)}{f(a)} \IE^{a} [f
(Z_{t_n-t_{l+1}})].\nonumber
\end{eqnarray}
In the first part of the proof, we provide an upper bound for
\[
\psi(a):=|(P_{l,l+1}-P_{t_{l+1}-t_l}) P_{t_n-t_{l+1}} f(a)|\qquad
\mbox{for $0\leq a< l$.}
\]
We couple $Z_{\Delta t_{l}}$ and $\cZ[l,l+1]$ under $\IP^a$ as in
Lemma~\ref{le0525-1} and denote by $\Upsilon^{\ssup1}$ and $\Upsilon
^{\ssup2}$
the respective random variables. There are two possibilities for the
coupling to fail: either $\Upsilon^{\ssup1}\geq a+2$ and $\Upsilon
^{\ssup2}=a+1$,
or $\Upsilon^{\ssup1} =a$ and $\Upsilon^{\ssup2}=a+1$. Consequently,
%
%
\begin{eqnarray}\label{eq0528-3}
\psi(a)&\leq&\IP\bigl( \Upsilon^{\ssup1} =a, \Upsilon^{\ssup2} =a+1\bigr)
\vphi
(a, a+1) \nonumber\\[-8pt]\\[-8pt]
&&{} + \IE\bigl[ \ind_{\{\Upsilon^{\ssup1}\geq a+1\}} \vphi\bigl(a+1,
\Upsilon^{\ssup1}
\bigr)\bigr].\nonumber
\end{eqnarray}
Since, by Taylor's formula,
\[
\IP\bigl(\Upsilon^{\ssup1}=a, \Upsilon^{\ssup2} =a+1\bigr)
= e^{-f(a) \Delta
t_l}-\bigl(1-f(a) \Delta t_l\bigr)\leq\tfrac12 (f(a) \Delta t_l)^2,
\]
we get for the first term of \eqref{eq0528-3}, using (\ref{eq0506-1}),
%
%
\begin{eqnarray}\label{eq0528-1}
&&\IP\bigl( \Upsilon^{\ssup1} =a, \Upsilon^{\ssup2} =a+1\bigr) \vphi(a,
a+1)\nonumber\\
&&\qquad\leq\frac12 (f(a) \Delta t_l)^2 \frac{\Delta f(a)}{f(a)} \IE
^{a} [f
(Z_{t_n-t_{l+1}})]\\
&&\qquad\leq f(a) (\Delta t_l)^2 \IE^{a} [f (Z_{t_n-t_{l+1}})].\nonumber
\end{eqnarray}
Now consider the second term in \eqref{eq0528-3}. We have
%
%
\begin{equation}\label{eq0528-2}\qquad
\IE\bigl[\ind_{\{\Upsilon^{\ssup1}\geq a+1\}}
\vphi\bigl(a+1,\Upsilon^{\ssup1}\bigr)\bigr] \leq
\underbrace{\IP\bigl(\Upsilon^{\ssup2}=a+1\bigr)}_{\leq f(a) \Delta t_l}
\IE^{a+1}[\vphi(a+1,Z_{\Delta t_l})].
\end{equation}
By Lemma~\ref{cor0602-1} we have $\IE^{a+1}[ f(Z_{\Delta t_l})]\leq
f(a+1) e^{\Delta f(a+1) \Delta t_l}$,
so that we conclude with (\ref{eq0506-1}) that
\begin{eqnarray*}
\IE^{a+1}[ \vphi(a+1,Z_{\Delta t_l})] &\leq& \bigl(e^{\Delta f(a+1) \Delta
t_l} -1\bigr) \IE^{a+1} [f(Z_{t_n-t_{l+1}}) ] \\
&\leq& 2 \Delta t_l \IE^{a+1} [f(Z_{t_n-t_{l+1}}) ],
\end{eqnarray*}
where we used in the last step that $\Delta f(a+1)<1$ and that $e^x\leq
1+2 x$ for $x\in[0,1]$.
We combine this with the
estimates (\ref{eq0528-3}), (\ref{eq0528-1}) and (\ref{eq0528-2}), and get
\[
\psi(a)\leq3 f(a) (\Delta t_l )^2 \IE^{a+1} [f(Z_{t_n-t_{l+1}}) ].
\]
In the next step, we deduce an estimate for $|\Sigma_l|$ defined in
\eqref{eq0423-2}. One has
\begin{eqnarray*}
|\Sigma_l|&\leq& P_{m,l} \psi(k)\leq3 \Delta t_{l} \IE^k\bigl[
\Delta t_l f(\cZ[m,l]) \IE^{\cZ[m,l]+1}[f(Z_{t_n-t_{l+1}})]\bigr]\\
&=& 3\Delta t_l \IE^k\bigl[\ind_{\{\Delta\cZ[m,l]=1\}} \IE^{\cZ
[m,l+1]}[f(Z_{t_n-t_{l+1}})]\bigr].
\end{eqnarray*}
By Lemma~\ref{ledom2} we get
%
%
\begin{eqnarray}\label{eq0423-1}
|\Sigma_l|&\leq&3\Delta t_l \IP^k(\Delta\cZ[m,l]=1) \IE^{k+1}
\bigl[\IE^{\cZ[m,l+1]}[f(Z_{t_n-t_{l+1}})]\bigr]\nonumber\\
&=& 3(\Delta t_l)^2 \IE^k[f(\cZ[m,l])] \IE^{k+1}\bigl[\IE^{\cZ
[m,l+1]}[f(Z_{t_n-t_{l+1}})]\bigr]\\
&=& 3(\Delta t_l)^2 P_{m,l} f(k) P_{m,l+1} P_{t_n-t_{l+1}}
f(k+1).\nonumber
\end{eqnarray}
We write $P_{t_n-t_{l+1}} f(k+1)=P_{t_{l+2}-t_{l+1}} P_{t_n-t_{l+2}}
f(k+1)$ and note that, by Lem\-ma~\ref{le0528-3},
$P_{t_n-t_{l+2}}f$ is concave. Therefore, we get with Lemma \ref
{le0528-2} that
$P_{t_n-t_{l+1}} f(k+1)\leq\tilde P_{l+1,l+2} P_{t_n-t_{l+2}}f(k+1)$.
Successive applications of this estimate and Lemma~\ref{le0528-1} yield
%
%
\begin{equation}\label{eq0529-1}
P_{m,l+1} P_{t_n-t_{l+1}} f(k+1)\leq\tilde P_{m,n}f(k+1)\leq
C_{\mmbox{\ref{le0528-1}}} P_{m,n} f(k+1).
\end{equation}
%
Recall from Lemma~\ref{le0602-1} that
$P_{m,l} f(k)\leq(\frac lm)^{\gamma^+} f(k)$. 
Combining with (\ref{eq0423-2}), (\ref{eq0423-1}) and (\ref{eq0529-1}) yields
%
%
\begin{eqnarray}\label{eq0529-2}
&&|P_{m,n} f(k)-P_{t_n-t_m} f(k)|\nonumber\\
&&\qquad\leq3 C_{\mmbox{\ref{le0528-1}}} f(k) P_{m,n} f(k+1) m^{-\gamma
^+} \sum
_{l=m}^{n-1} l^{-2+\gamma^+}\\
&&\qquad\leq C_{\mmbox{\ref{SGapprox}}} \frac{f(k)}m P_{m,n}
f(k),\nonumber
\end{eqnarray}
for a suitably defined constant $C_{\mmbox{\ref{SGapprox}}}$ depending
only on
$f$, as required. 
\end{pf*}

\section{\texorpdfstring{Properties of the family $(A_\alpha\dvtx0<\alpha<1)$ of operators}
{Properties of the family (A alpha: 0<alpha<1) of operators}}\label{se3}

The objective of this section is to study the operators $A_\alpha$ and
relate them to
the tree $\mathrm{INT}$. We start with two lemmas on the functional
analytic nature
of the family $(A_\alpha\dvtx\alpha\in{\mathcal I})$.
\begin{lemma}\label{propop}
\textup{(a)} For any $0<\alpha<1$ the following are equivalent:
\begin{longlist}[(ii)]
\item[\hphantom{(bi)}(i)] $A_\alpha1(0)<\infty$;
\item[\hphantom{(b)}(ii)] $A_\alpha g \in{\mathbf C}({\mathcal S})$ for all $g\in
{\mathbf C}({\mathcal S})$.
\end{longlist}

The set of $\alpha$ where these conditions hold is denoted
by ${\mathcal I}$.

\textup{(b)} For any $\alpha\in{\mathcal I}$ the operator $A_\alpha$ is
strongly positive.

\textup{(c)} For any $\alpha\in{\mathcal I}$ the operator $A_\alpha$ is
compact.
\end{lemma}
\begin{pf}
Recalling the Arzel\`a--Ascoli theorem,
the only nontrivial claim is that, if $A_\alpha1(0)<\infty$, then the
family $(A_\alpha g \dvtx\|g\|_\infty<1)$ is
equicontinuous. To this end recall that, for $\tau\leq\sigma\leq\infty
$, by Remark~\ref{redom}, we have $\M^\tau\geq M^\sigma$ and hence
\[
| A_\alpha g(\tau) - A_\alpha g(\sigma) |
\leq\int_0^\infty e^{-\alpha t} d(\M^{\tau}-\M^{\sigma})(t) .
\]
Equicontinuity at $\infty$ follows from this by recalling the
definition $\M^\infty=\lim_{\tau\uparrow\infty} \M^\tau$.
Elsewhere, for $\sigma<\infty$, we use the straightforward coupling of
the processes $(Z^{\tsupp\tau}_t\dvtx t\geq0)$ and
$(Z^{\tsupp\sigma}_t\dvtx t\geq0)$ with the property that \textit{if}
$Z_{\sigma-\tau}^{\tsupp\sigma}=0$ then
$Z_t^{\tsupp\tau}=Z_{t+{\sigma-\tau}}^{\tsupp\sigma}$.

Hence we get
%
%
\begin{eqnarray}\label{eq1008-5}
\int_0^\infty e^{-\alpha t} \,d(\M^{\tau}-\M^{\sigma})(t) &\leq&
\bigl(1-e^{-\alpha(\sigma-\tau)}\bigr) \int_0^\infty e^{-\alpha t} \,d\M^\tau
(t) \nonumber\\[-8pt]\\[-8pt]
&&{} + \IE\biggl[ \int_0^\infty e^{-\alpha t} dZ^{\tsupp\tau}_t \ind\bigl\{
Z_{\sigma-\tau}^{\tsupp\sigma}>0 \bigr\}\biggr]
.\nonumber
\end{eqnarray}
Note that $\int_0^\infty e^{-\alpha t} \,d\M^\tau(t)\leq\IE[ \int
_0^\infty e^{-\alpha t} dZ^{\tsupp\tau}(t)] \leq A_\alpha1 (0)<\infty$,
and that\break $\IP(Z_{\sigma-\tau}^{\tsupp\sigma}>0)\leq{\IP^1(Z_{\sigma-\tau
}>1)}\downarrow0$ as $\sigma\downarrow\tau$. Hence, both terms on the
right-hand side of (\ref{eq1008-5})
can be made small by making $\sigma-\tau$ small, proving the claim.
\end{pf}
\begin{lemma}\label{evconvex}
The function $\alpha\mapsto\log\rho(A_\alpha)$ is convex on $\mathcal I$.
\end{lemma}
\begin{pf}
By Theorem 2.5 of~\cite{Ka82} the function $\alpha\mapsto\log\rho
(A_\alpha)$ is convex, if for each
positive $g\in{\mathbf C}({\mathcal S})$, $\eps>0$ and triplet $\alpha
_1\leq\alpha_0 \leq\alpha_2$ in ${\mathcal I}$,
there are finitely many positive $g_j\in{\mathbf C}({\mathcal S})$ and
functions $\phi_j\dvtx\mathcal I \to\IR$, $j\in\{1,\ldots,m\}$,
with $\log\phi_j$ convex, such that
\[
\Biggl\| A_{\alpha_k}g-\sum_{j=1}^m \phi_j(\alpha_k)g_j \Biggr\| \leq\eps
\qquad\mbox{for all }k\in\{0,1,2\}.
\]
This criterion is easily checked using the explicit form of $A_\alpha$,
$0<\alpha<1$.
\end{pf}

With the help of the following lemma, Theorem~\ref{main1} follows from
Theorem~\ref{main2}. The result is a variant of a standard result in
the theory of
branching random walks adapted to our purpose; see, for example, Hardy
and Harris~\cite{HH09} for a good account of the general theory.
\begin{lemma}\label{opimp}
The $\mathrm{INT}$ dies out almost surely if and only if there exists
$0<\alpha<1$
such that $A_\alpha$ is a compact linear operator with spectral radius
$\rho(A_\alpha)\leq1$.
\end{lemma}
\begin{pf} 
Suppose that such an $\alpha$ exists. By the Krein--Rutman theorem
(see, e.g., Theorem 1.3 in Section 3.2 of~\cite{Pi95})
there exists a eigenvector $v\dvtx{\mathcal S} \to[0,\infty)$
corresponding to the eigenvalue $\rho(A_\alpha)$.
Our operator $A_\alpha$ is strongly positive, that is, for every $g\geq
0$ which is positive somewhere,
we have
\[
\min_{\tau\in{\mathcal S}} A_\alpha g(\tau)>0,
\]
so that $v$ is also bounded away from zero. Let $Y^{\ssup
n}_{\tau}(dt\,
dx)$ be the empirical measure of types and
positions of all the offspring in the $n$th generation of an IBRW
started by a single particle of type $\tau$ positioned at the origin.
With every generation of particles in the IBRW we associate a score
\[
X_n:=\int Y^{\ssup n}_{\tau}(dt \, dx) e^{-\alpha x} \frac
{v(t)}{v(\tau)}.
\]
%
The assumption $\rho(A_\alpha)\leq1$
implies that $(X_n \dvtx n\in\IN)$ is a supermartingale and thus
almost surely convergent.
Now fix some $N>1$, an integer $n\geq2$ and the state at
generation $n-1$. Suppose there
is a particle with location $x<N$ in the $(n-1)$st generation. Then
there is a positive probability
(depending on $N$ but not on $n$) that $X_n-X_{n-1}>1$ and, as $(X_n
\dvtx n\in\IN)$
converges, this can only happen for finitely many $n$. Hence the location
of the leftmost particle in the IBRW diverges to $+\infty$ almost
surely. This implies that the INT dies
out almost surely.

Conversely, we assume that $\mathcal I$ is nonempty and fix $\alpha\in
{\mathcal I}$.
The Krein--Rutman theorem gives the existence of an eigenvector\vadjust{\goodbreak} of the
dual operator,
which is a positive, finite measure $\nu$ on the type space ${\mathcal
S}$ such that
$\int v(t) \nu(dt)=1$ and, for all continuous, bounded $f\dvtx
{\mathcal S} \to\IR$,
\[
\int A_\alpha f(t) \nu(dt) = \rho(A_\alpha) \int f(t) \nu
(dt) .
\]
Because $A_\alpha$ is a strongly positive operator, the Krein--Rutman
theorem implies that there exists
$\lambda_0<\rho(A_\alpha)$ such that $|\lambda|\leq\lambda_0$ for all
$\lambda\in\sigma(A_\alpha)\setminus
\{\rho(A_\alpha)\}$, where $\sigma(A_\alpha)$ denotes the spectrum of
the operator. Hence $\rho(A_\alpha)$
is separated from the rest of the spectrum and by Theorem IV.3.16
in~\cite{Kato} this holds for all parameters
in a small neighborhood of $\alpha$. Hence, arguing as in Note 3 on
Chapter II in~\cite{Kato}, pages 568 and 569, the mapping
$\alpha\mapsto\rho(A_\alpha)$ is differentiable and its derivative equals
%
%
\begin{equation}\label{deriv}
\rho'(A_\alpha):=\frac{d}{d\alpha} \int A_\alpha v(t) \nu(dt)
= \int\frac{\partial}{\partial\alpha} A_\alpha v(t) \nu(dt),
\end{equation}
where the second equality can be inferred from the minimax
characterisation of eigenvalues;
see, for example, Theorem 1 in~\cite{Ra83}. Given $\tau\in{\mathcal S}$
we define a martingale by
\[
W^{\ssup n}_\tau= \rho(A_\alpha)^{-n} \iint\frac{v(t)}{v(\tau)}
e^{-\alpha x} Y^{\ssup n}_{\tau}(dt \,dx)
\]
and
argue as in Theorem 1 of~\cite{KRS} that it converges almost surely to
a strictly positive limit $W_\tau$ if
%
%
\begin{equation}\label{marco}
\log\rho(A_\alpha) - \frac{\alpha \rho'(A_\alpha)}{\rho(A_\alpha)} >0
\quad\mbox{and}\quad \sup_{\tau\in\mathcal S}\IE\bigl[ W^{\ssup
1}_\tau\log W^{\ssup1}_\tau\bigr]<\infty.
\end{equation}
Let us assume for the moment that the second condition holds true for
all $\alpha\in\mathcal I$. Then,
if $\alpha$ is such that the limit $W_\tau$ exists and is positive, it
also exists for the offspring of
any particle of type $\tau$ in position $x$, and we denote it
by $W_\tau(x)$. By decomposing the population in the $m$th generation
according to their ancestor in the $n$th generation,
and then letting $m\to\infty$, we get
\[
W_\tau= \rho(A_\alpha)^{-n} \int\frac{v(t)}{v(\tau)} e^{-\alpha
x} W_t(x) Y^{\ssup n}_{\tau}(dt \, dx) .
\]
Denoting by $P_\tau$ the law of the IBRW started with a particle at the
origin of type $\tau$, we now look at the
IBRW under the changed measure
\[
dQ = \int\nu(d\tau) v(\tau) W_\tau \,dP_\tau .
\]
Given a sample IBRW we build a measure $\mu$ on the set of all infinite
sequences
\[
((x_0,t_0), (x_1,t_1), \ldots),
\]
where $x_j$ is the location and $t_j$ the type of a particle in the
$j$th generation,
which is a child of a particle in position $x_{j-1}$ of type $t_{j-1}$,
for all $j\geq1$. This\vadjust{\goodbreak} measure is determined by the
requirement that, for any permissible sequence
\begin{eqnarray*}
&&\mu\{ ((y_0,s_0), (y_1,s_1),\ldots) \dvtx y_0=x_0, s_0=t_0,
\ldots, y_n=x_n, s_n=t_n\} \\
&&\qquad= \rho(A_\alpha)^{-n} \frac{v(t_n)}{v(t_0)} \exp\{-\alpha (x_n-x_0)
\}
\frac{W_{t_n}(x_n)}{W_{t_0}(x_0)}. 
\end{eqnarray*}
%
Looking unconditionally at the random sequence of particle types thus
generated, we note that it is a stationary
Markov chain on ${\mathcal S}$ with invariant distribution $v(t) \nu
(dt)$ and transition kernel given by
\begin{eqnarray*}
P_{t_0}(\ell)& = &\rho(A_\alpha)^{-1} \frac{v(\ell)}{v(t_0)} \int
_0^\infty e^{-\alpha t} \,d\M^{t_0}(t) ,\\
P_{t_0}(dt)& = &\rho(A_\alpha)^{-1} \frac{v(t)}{v(t_0)} e^{\alpha
t} \,d\M(t) \qquad\mbox{for }t\geq0 .
\end{eqnarray*}
Using first Birkhoff's ergodic theorem and then \eqref{deriv} we see
that, $Q$-almost surely,
$\mu$-almost every path has speed 
\begin{eqnarray*}
\lim_{n\to\infty} \frac{x_n}{n}
&=& \frac1{\rho(A_\alpha)} \int\IE\biggl[ \int Y^{\ssup1}_{t_0}(dt\,
dx) x e^{-\alpha x} \frac{v(t)}{v(t_0)} \biggr] v(t_0) \nu
(dt_0) \\
&=& - \frac1{\rho(A_\alpha)} \int\frac{\partial}{\partial
\alpha}\, \frac{A_\alpha v(t_0)}{v(t_0)} v(t_0) \nu(dt_0)
\\
&=& - \frac{\rho'(A_\alpha)}{\rho(A_\alpha)} = -\frac{d}{d\alpha} \log
\rho(A_\alpha).
\end{eqnarray*}
Suppose that 
$\alpha_0\in{\mathcal I}$ is such that
\[
\rho(A_{\alpha_0}) =\min_{\alpha\in{\mathcal I}} \rho(A_\alpha) >1.
\]
From Lemma~\ref{evconvex} we can infer that there exists $\alpha>\alpha
_0$ such that the first condition in \eqref{marco} holds and
\[
-\frac{d}{d\alpha} \log\rho(A_\alpha)<0.
\]
This implies that, $Q$-almost surely, there exists an ancestral line of
particles diverging to $-\infty$. For the
IBRW started with a particle at the origin of type $\ell$, we therefore
have a positive probability that an ancestral
line goes to $-\infty$. This implies that the INT has a positive
probability of survival.

To ensure that the second condition in \eqref{marco} holds, we can use
a cut-off procedure and replace the
offspring distribution $Y^{\ssup1}(dt \, dx)$ by one that takes only
the first $N$ children to the right
and left into account. It is easy to see that, for fixed $0<\alpha<1$
and sufficiently large $N$, we can ensure
that the modified operator $A_\alpha^{\ssup N}$ is close to the
original one in the operator norm, and as large as we wish if the
original operator is ill defined. Hence the continuity of the spectral
radius in the operator norm
ensures that $\lim_{N\to\infty}\rho(A^{\ssup N}_\alpha
)=\rho(A_\alpha)$, with\vadjust{\goodbreak} the spectral radius of an ill-defined operator
being infinity. Using Lemma~\ref{evconvex} and the fact that a sequence
of convex functions, which converges pointwise, converges uniformly on
every closed set, we
can choose $N$ so that for all $0<\alpha<1$ the modified operators satisfy
$\rho(A^{\ssup N}_\alpha)>1$, while the cut-off ensures that the second
criterion in \eqref{marco} automatically holds.
The argument above can now be applied and yields the existence of an
ancestral line of particles diverging to $-\infty$,
which then automatically also exists in the original IBRW.\vspace*{-2pt}
\end{pf}

Our proofs, in particular the crucial sprinkling technique, rely on the
following continuity property of the survival
probability $p(f)$ of the INT for the attachment rule $f$.\vspace*{-2pt}
%
\begin{lemma}\label{pcontinuous}
One has
\[
\lim_{\eps\dto0} p(f-\eps)= p(f).\vspace*{-2pt}
\]
\end{lemma}
\begin{pf}
We only need to consider the case where $p(f)>0$, as otherwise both
sides of the equation are zero. We denote
by $\rho(\alpha, f)$ the spectral radius of the operator $A_\alpha$
formed with respect to the attachment function $f$,
setting it equal to infinity if the operator is ill defined.
The assumption $p(f)>0$ implies, by Lem\-ma~\ref{opimp}, that for all
$0<\alpha<1$ we have $\rho(\alpha, f)>1$.
As the operator norm $\|A_\alpha\|$ for the operator formed with
respect to the attachment function $f-\eps$ depends
continuously on $\eps\geq0$, we can use the continuous dependence of
the spectral radius on the operator norm to obtain,
for fixed~$\alpha$,
\[
\lim_{\eps\downarrow0} \rho(\alpha,f-\eps)=\rho(\alpha,f).
\]
As a sequence of convex functions, which converges pointwise, converges
uniformly on every closed set, we
find $\eps>0$ such that $\rho(\alpha, f-\eps)>1$ for all $0<\alpha
<1$. Thus,
using again Lemma~\ref{opimp}, we have $p(f-\eps)>0$.

Now we look at the IBRW started with one particle of type $\ell$ in
position $t$, constructed using
the attachment rule $f-\eps$, such that any particle with position $>0$
is killed along with its offspring.
We denote by $E(\eps,t)$ the event this process survives forever, and
by $V(\eps,t, \kappa)$ the event
that a particle reaches a site $<\kappa$. Then we have
\[
\lim_{\kappa\to-\infty} \inf_{t\le\kappa} \IP(E(\eps,t))=1.
\]
For fixed $\kappa<0$ and $0\leq\eps\leq\eps_0$ we have
\[
\IP( E(\eps,t) ) \geq\IP( V(\eps,t,\kappa) ) \IP
( E(\eps_0,\kappa) )
\stackrel{\eps\downarrow0}{\longrightarrow} \IP( V(0,t,\kappa)
) \IP( E(\eps_0,\kappa) ).
\]
Note that the first probability on the right is greater or equal to
$p(f)$ and that the second probability
tends to one, as $\kappa$ tends to $-\infty$.\vspace*{-2pt}
\end{pf}


\section{The giant component}\label{se4}

This section provides two crucial tools: a tool to obtain global
results from our local approximations of neighborhoods given by
the ``sprinkling'' argument in\vadjust{\goodbreak} Proposition~\ref{sprinkling}, and an a
priori lower bound on the size of the connected components of
the oldest vertices in the system given in Lem\-ma~\ref{initial}. We
follow the convention that a sequence of events depending on
the index $N$ holds \textit{with high probability} if the probability of
these events goes to one as $N\uparrow\infty$.
\begin{prop}[(Sprinkling argument)]\label{sprinkling}
Let $\eps\in(0,f(0))$, $\kappa>0$, and $\bar f(k)= f(k)-\eps$ for
integers $k\geq0$.
Suppose that $(c_N)_{N\in\IN}$ is a sequence of integers with
\[
\lim_{N\uparrow\infty} \biggl[\sfrac{1}{2} \kappa\eps c_N-\log N
\biggr]=\infty\quad\mbox{and}\quad \lim_{N\to\infty} \frac{c_N^2}{N}=0,
\]
and that, for the preferential attachment graphs $(\bar\cG_N)_{N\in\IN
}$ with attachment rule $\bar f$, we have
\[
\sum_{v=1}^N \ind{\{|\bar{\mathcal C}_N(v)|\geq2c_N\}}\geq\kappa N
\qquad\mbox{with high probability,}
\]
where $\bar{\mathcal C}_N(v)$ denotes the connected component of the
vertex $v$ in $\bar\cG_N$.
Then there exists a coupling of the graph sequences $(\cG_N)_{N\in\IN}$
with $(\bar\cG_N)_{N\in\IN}$ such that $\bar\cG_N\leq\cG_N$
and all connected components of $\bar\cG_N$ with at least $2c_N$ vertices
belong to one connected component in $\cG_N$ with at least $\kappa N$
vertices, with high probability.
\end{prop}
\begin{pf}
Note that we can couple $\bar\cG_N$ and an independent Erd\H{o}s--R\'
enyi graph 
$\cG_N^{\ER}$ with edge probability $\eps/N$ with $\cG_N$ such that
%
%
\begin{equation}\label{eq0327-1}
\bar\cG_N\leq\bar\cG_N \vee\cG_N^{\ER} \leq\cG_N.
\end{equation}
Here, $\bar\cG_N\vee\cG_N^{\ER}$ denotes the graph in which all edges
are open that are open in at least one of the two graphs,
and $\cG'\leq\cG''$ means that all edges\vspace*{2pt} that are open in $\cG'$ are
also open in $\cG''$.
We denote by $V'_N$ the vertices in $\bar\cG_N$ that belong to
components of size at least $2c_N$
and write $V'_N$ as the disjoint union $C_1\cup\cdots\cup C_M$, where
$C_1,\ldots,C_M$ are sets of vertices
such that:
\begin{itemize}
\item$|C_j|\in[c_N,2c_N]$ and
\item$C_j$ belongs to one component in $\bar\cG_N$, for each
$j=1,\ldots,M$.
\end{itemize}
Recall (\ref{eq0327-1}), and note that given $\bar\cG_N$ and the sets
$C_1,\ldots, C_M$,
the Erd\H{o}s--R\'enyi graph $\cG_N^{\ER}$ connects two distinct sets
$C_i$ and $C_j$ with probability at least
\[
p_N:=1-\biggl(1-\sfrac{\eps}{N}\biggr)^{c_N^2}\geq1- e^{-(\eps/N)
c_N^2}\sim\frac\eps N c_N^2.
\]
By identifying the individual sets as one vertex and interpreting the
$\cG_N^{\ER}$-connections as edges, we obtain a new random graph.
Certainly, this dominates an Erd\H{o}s--R\'enyi graph with $M$ vertices
and success probability $p_N$, which has edge intensity\vadjust{\goodbreak} $M p_N$.
By assumption,\vspace*{2pt} $\frac12 \frac{\kappa N}{c_N}\leq M\leq N$ with high
probability. Hence
$M\to\infty$ and $M p_N-\log M\to\infty$ in probability as $N\uparrow
\infty$.
By~\cite{Ho09}, Theorem 5.6, the new Erd\H{o}s--R\'enyi graph is
connected with high probability.
Hence, all vertices of $V_N'$ belong to one connected component in $\cG
_N$, with high probability.
\end{pf}

We need an ``a priori'' argument asserting that the connected components
of the old vertices are
large with high probability. This will, in particular, ensure that the
connected component of any vertex connected
to an old vertex is large.
\begin{lemma}[(A priori estimate)]\label{initial}
Let $(c_N)_{N\in\IN}$ and $({n_N})_{N\in\IN}$ be sequences of positive
integers such that
\[
\lim_{N\to\infty} \frac{c_N}{\log N\log\log N} =0 \quad\mbox{and}\quad \lim
_{N\to\infty} \frac{\log{n_N}}{\log N}=0.
\]
Denote by ${\mathcal C}_N(v)\subset{\mathcal G}_N$ the connected
component containing $v\in\{1,\ldots,N\}$. Then
\[
\IP\bigl(\#{\mathcal C}_N(v)< c_N \mbox{ for any } v\in\{1,\ldots
,{n_N}\} \bigr) \longrightarrow0.
\]
\end{lemma}
\begin{pf} 
We only need to show this for the case when $f$ is constant, say equal
to $\beta>0$, as all
other cases stochastically dominate this one. Note that in this case
\textit{all} edge probabilities are
independent. We first fix a vertex $v \in\{1,\ldots
,{n_N}\}$ and denote by $X_1=X_1(v)$ the number of its direct neighbours
in $({n_N},N/\log N]$. We obtain, for any $\lambda>0$,
\[
\IE e^{-\lambda X_1} = \prod_{j={n_N}}^{\lfloor N/\log N \rfloor-1}
\biggl( \sfrac{\beta}{j} e^{-\lambda} +
\biggl(1-\sfrac{\beta}{j}\biggr) \biggr),
\]
and hence, for sufficiently large $N$,
\[
\log\IE e^{-\lambda X_1}
\leq-\beta (1-e^{-\lambda}) \sum_{j={n_N}}^{\lfloor N/\log
N \rfloor-1} \frac1j
\leq- \sfrac{3}{4} \beta (1-e^{-\lambda}) \log N .
\]
By the exponential Chebyshev inequality we thus get for sufficiently
large~$N$,
%
%
\begin{equation}\label{Nesti}
\IP\biggl( X_1 < \sfrac{\beta}{2} \log N\biggr) \leq
N^{\lambda{\beta}/2 - ({3\beta}/4) (1-e^{-\lambda})} \leq
N^{-{\beta}/{32}}
\end{equation}
choosing $\lambda=\frac12$ and using that $1-e^{-x} \geq x-\frac12
x^2$ for $x\geq0$ in the last step. 
Now
let $X_2=X_2(v)$ be the number of direct neighbors in $(N/\log N,N]$ of
any of the $X_1(v)$ vertices who are
direct neighbors of $v$ in $({n_N},N/\log N]$.
Since by assumption $f(k)=\beta$ for all $k$, we obtain, for any
$\lambda>0$,
\[
\IE[ e^{-\lambda X_2} | X_1 ]
= \prod_{j=\lfloor N/\log N\rfloor}^{N-1} \biggl( 1+ (e^{-\lambda}-1)
\biggl(1-\biggl(1-\sfrac{\beta}{j}\biggr)^{X_1}\biggr)\biggr)\vadjust{\goodbreak}
\]
and hence, for sufficiently large $N$, on the event $\{X_1\geq\sfrac
{\beta}{2} \log N\}$,
\begin{eqnarray*}
\log\IE[ e^{-\lambda X_2} | X_1 ]
& \leq&-(1-e^{-\lambda}) \sfrac{3\beta}{4} X_1 \sum_{j=\lfloor
N/\log N\rfloor}^{N-1} \frac1j
\\
&\leq&-(1-e^{-\lambda}) \sfrac{\beta^2}{4} \log N \log\log N.
\end{eqnarray*}
By \eqref{Nesti} and the exponential Chebyshev inequality (with $\lambda
=1$) we thus get for sufficiently large $N$,
\begin{eqnarray*}
\IP\bigl(X_2(v) < c_N \bigr)
&\leq&\IP\biggl( X_1 < \sfrac{\beta}{2} \log N\biggr)
+ \IP\biggl( X_2(v) < c_N \Big| X_1\geq\sfrac{\beta}{2} \log N \biggr)\\
&\leq& N^{-{\beta}/{32}} + N^{- ({\beta^2}/8)\log\log N +
c_N/\log N}.
\end{eqnarray*}
Let $\lambda=\frac12$. By our assumptions on $(c_N)_{N\in\IN}$ and
$({n_N})_{N\in\IN}$ the sum of the
right-hand sides over all $v \in\{1,\ldots,{n_N}\}$ goes to zero,
ensuring that
$\#{\mathcal C}_N(v) \geq X_2(v) \geq c_N$ for all $v \in\{1,\ldots
,{n_N}\}$ with
high probability.
\end{pf}

\section{The exploration process}\label{se1stcoupling}

Our aim is to ``couple'' certain aspects of the network to an easier
object, namely a random 
tree. To each of these objects we associate a dynamic process called
the {exploration process}.
In general, an \textit{exploration process} of a graph successively
collects information about the connected component
of a fixed vertex by following edges emanating from already discovered
vertices in a well-defined order, so that at each
instance the explored part of the graph is a connected subgraph of the cluster.
We show that the exploration processes of the network and the labeled tree
can be defined on the same probability space in such a way that up to a
stopping time, which is typically large,
the explored part of the network and the tree coincide.

\subsection{A random labeled tree}\label{se51}

We now describe a tree ${\mathbb T}(w)$ which informally describes the
neighborhood of a vertex $w\in\cG_N$.
Any vertex in the tree is labeled by two parameters: its \textit
{location}, an element of $\{1,\ldots,N\}$, and its \textit{type},
an element of $\{\ell\}\cup\{1,\ldots,N\}$. The root is given as a
vertex with location $w$ and type $\ell$. A vertex $v$ with
location $i$ and type $\ell$ produces independently descendants in the
locations $1,\ldots, i-1$ (i.e., to its left) of type
$i$ with probability
\[
\IP(v \mbox{ has a descendant in }j \mbox{ of type }i)=\IP(\Delta\cZ[j,i-1]=1).
\]
Moreover, independently it produces descendants to its right, which are
all of type~$\ell$, in such a way that
the cumulative sum of these descendants is distributed according to the
law of $(\cZ[i,j] \dvtx i+1\leq j\leq n)$.
A vertex $v$ of type $k$ produces descendants to the left in the same
way as a vertex of type $\ell$, and independently
it produces descendants to the right, which are all of type $\ell$, in
such a way that the cumulative sum of these descendants
is distributed as $( \cZ[i,j]-\ind_{[k,\infty)}(j) \dvtx i+1\leq j\leq
n)$ conditioned on $\Delta\cZ[i,k-1]=1$.

Observe that, given the tree and the locations of the vertices, we may
reconstruct the types of the vertices in a deterministic
way: any vertex whose parent is located to its left has the type $\ell
$, otherwise the type of the vertex is the location
of the parent.

The link between this labeled tree and our network is given in the
following proposition, which will be proved in Section~\ref{se53}.
\begin{prop}\label{corcoup1}
Suppose that $(c_N)_{N\in\IN}$ is a sequence of integers with
\[
\lim_{N\to\infty} \frac{c_N}{\log N \log\log N}=0.
\]
Then one can couple the pair $(V,\cG_N)$ consisting of the network and
a uniformly chosen vertex $V$ with
$\IT(V)$ such that with high probability
\[
\#\cC_N(V) \wedge c_N= \#\IT(V)\wedge c_N.
\]
%
\end{prop}

\subsection{Exploration of the network}

We now specify how we explore a graph like our network or the tree
described above, that is, we specify
the way we collect information about the connected component, or
cluster, of a particular vertex $v$.
In the first step, we explore all immediate neighbors of $v$ in the graph.
To explain a general exploration step we classify the vertices in three
categories:
\begin{itemize}
\item\textit{veiled vertices}: vertices for which we have not yet found
connections to the cluster of $v$;
\item\textit{active vertices}: vertices for which we already know that
they belong to the cluster, but for which we have not yet explored all
its immediate neighbors;
\item\textit{dead vertices}: vertices which belong to the cluster and
for which all immediate neighbors have been explored.
\end{itemize}

After the first exploration step the vertex $v$ is marked as dead, its
immediate neighbors as active and all
the remaining vertices as veiled.
In a general exploration step, we choose the \textit{leftmost} active
vertex, set its state to \textit{dead}, and explore its immediate
neighbors. The newly found \textit{veiled} vertices are marked as
\textit
{active}, and we proceed with another exploration step until there are
no active vertices left.

In the following, we couple the exploration processes of the network
and the random labeled tree started with a particle at position $v$ and
type $\ell$ up to a stopping time $T$. Before we introduce the coupling
explicitly, let us quote adverse events which stop the coupling.
Whenever the exploration process of the network revisits an active
vertex we have found a cycle 
in the network. We call this event (E1) and
stop the exploration so that, before time $T$, the explored part of the
neighborhood of $v$ is a tree with each node having a unique location.
Additionally, we stop once the explored part of the network differs
from the explored part of the random labeled tree, calling this event (E2),
we shall see in Section~\ref{se53} how this can happen. In cases (E1)
and (E2) we say that \textit{the coupling fails}.


Further reasons to stop the exploration are, for parameters $n_N,c_N\in
\IN$ with
$1\leq n_N,c_N\leq N$:
\begin{longlist}[(C)]
\item[(A)] the number of dead and active vertices exceeds $c_N$;
\item[(B)] one vertex in $\{1,\ldots,n_N\}$ is activated;
\item[(C)] there are no more active vertices left.
\end{longlist}
If we stop the exploration without (E1) and (E2) being the case, we say
that the \textit{coupling succeeds}.
Once the exploration has stopped, the veiled parts of the random tree
and the network may be generated
independently of each other with the appropriate probabilities. Hence,
if we succeed in coupling the explorations,
we have coupled the random labeled tree and the network.


\subsection{Coupling the explorations}
\label{se53}

To distinguish both exploration processes, we use the term \textit
{descendant} for a child in the labeled random tree and
the term \textit{immediate neighbor} in the context of the neighborhood
exploration in the network.
In the initial step, we explore all immediate neighbors of $v$ and all
the descendants of the root. Both explorations are
identically distributed and they therefore can be perfectly coupled.
Suppose now that we have performed $k$ steps and that we have not yet
stopped the exploration. In particular, this means that both explored
subgraphs coincide and that any unveiled (i.e., active or dead) element
of the labeled random tree can be uniquely referred to by its location.
We now explore the descendants and immediate neighbors of the \textit
{leftmost} active vertex, say $n$.

\textit{First}, we explore the descendants to the left (veiled and dead)
and immediately check whether they themselves have right
descendants in the set of dead vertices. If we discover no dead
descendants, the set of newly found left descendants is
identically distributed to the immediate left neighbors in the network.
Thus we can couple both explorations such that
they agree in this case. Otherwise we stop the exploration due to (E2).

\textit{Second}, we explore the descendants to the right. If the
vertex $n$ is \textit{not of type $\ell$}, then we know already that $n$
has no \textit{right} descendants that were marked as dead as $n$ itself
was discovered. Since we always explore the leftmost active vertex
there are no new dead vertices to the right of $n$. Therefore, the
explorations to the right in the network and the random labeled tree
are identically distributed and we stop if we find right neighbors in
the set of active vertices due to\vadjust{\goodbreak} (E1). 
If the vertex $n$ is \textit{of type $\ell$}, then we have not gained any
information about its right \textit{descendants} yet.
If we find no right descendants in the set of dead vertices, it is
identically distributed to the immediate right neighbors of $n$ in the
network. We stop if right descendants are discovered that were marked
as dead, corresponding to (E2), or if right descendants are discovered
in the set of active vertices, corresponding to~(E1).
%
%
\begin{lemma}\label{premature}
Suppose that $(c_N)_{N\in\IN}$, $(n_N)_{N\in\IN}$ are sequences of
integers such that
\[
\lim_{N\to\infty} \frac{c_N^2}{n_N^{1-\gamma^+}}=0.
\]
Then the coupling of the exploration processes 
satisfies
\[
\lim_{N\to\infty} \sup_{v\in\{{n_N}+1,\ldots,N\}}
\IP\bigl(\mbox{coupling with initial vertex $v$ ends in $\mathrm{(E1)}$
or $\mathrm{(E2)}$}\bigr) =0,
\]
that is, the coupling succeeds with high probability.
\end{lemma}
\begin{pf} 
We analyze one exploration step in detail. Let $\mathfrak a$ and
$\mathfrak d$ denote the active and dead vertices of a feasible
configuration at the beginning of an exploration step, that is,
$\mathfrak a, \mathfrak d$ denote two disjoint subsets of $\{
{n_N}+1,\ldots,N\}$ with $\# (\mathfrak a\cup\mathfrak d)<c_N$ and
$\mathfrak a\not=\varnothing$.

The exploration of the minimal vertex $n$ in the set $\mathfrak a$ may
only fail 
for one of the following reasons:
\begin{longlist}[(Ib)]
\item[(Ia)] the vertex $n$ has left descendants in $\mathfrak d$,
\item[(Ib)] the vertex $n$ has left descendants which themselves have
right descendants in $\mathfrak d$ or
\item[(II)] the vertex $n$ has right descendants in $\mathfrak a\cup
\mathfrak d$.
\end{longlist}
Indeed, if (Ia) and (Ib) do not occur, then the exploration to the left
ends neither in state (E1) nor (E2),
and if (II) does not happen, the exploration to the right does not fail.

Conditionally on the configuration $(\mathfrak a,\mathfrak d)$, the
probability for the event (Ia) equals
%
\[
\IP( \exists{a\in\mathfrak d} \mbox{ such that } \Delta\cZ
[a,n-1]=1)\leq\sumtwo{a\in\mathfrak d}{a<n} \IP(\Delta\cZ[a,n-1]=1),
\]
whereas the probability for (Ib) is by Lemma~\ref{ledom2} equal to %
%
\begin{eqnarray*}
&&\IP( \exists a\in\mathfrak d^c\mbox{ and }b\in\mathfrak d\mbox{ such
that } \Delta\cZ[a,n-1]=\Delta\cZ[a,b-1]=1)\\
&&\qquad\leq\sumtwo{a\in\mathfrak d^c}{a<n} \sumtwo{b\in\mathfrak
d}{b>a} \IP(\Delta\cZ[a,n-1]=\Delta\cZ[a,b-1]=1)\\
&&\qquad\leq\sumtwo{a\in\mathfrak d^c}{a<n} \sumtwo{b\in\mathfrak d}{b>a}
\IP(\Delta\cZ[a,n-1]=1) \IP^1(\Delta\cZ[a,b-1]=1).
\end{eqnarray*}
If the vertex $n$ is of type $\tau\not=\ell$, then the conditional
probability of (II) is
\begin{eqnarray*}
&&\IP(\exists a\in\mathfrak a \mbox{ such that } \Delta\cZ
[n,a-1]=1
|\Delta\cZ[n,\tau-1]=1, \Delta\cZ[n,b-1]=0\\
&&\hspace*{271pt}\forall b\in\mathfrak
d\setminus\{\tau\})\\
&&\qquad\leq C_{\mmbox{\ref{condaway}}} \sumtwo{a\in\mathfrak a\cup
\mathfrak
d}{a>n} \IP^1( \Delta\cZ[n,a-1]=1)
\end{eqnarray*}
using first Lemma~\ref{condaway}
and then Lemma~\ref{ledom2}.

If the vertex $n$ is of type $\ell$, the conditional probability of
(II) is
%
\[
\IP(\exists a\in\mathfrak a\cup\mathfrak d \mbox{ such that }
\Delta\cZ[n,a-1]=1)
\leq\sumtwo{a\in\mathfrak a\cup\mathfrak d}{a>n} \IP( \Delta\cZ[n,a-1]=1).
\]
Since, by Lemma~\ref{le0416-1}, for any $a>n$,
\[
\IP^1(\Delta\cZ[n,a-1]=1)\leq\IP^1(\Delta\cZ
[{n_N}+1,{n_N}+1]=1),
\]
we conclude that the probabilities of the events (Ia) and (II) are
bounded by
\[
(2+C_{\mmbox{\ref{condaway}}}) c_N \IP^1(\Delta\cZ
[{n_N}+1,{n_N}+1]=1),
\]
independently of the type $\tau$. Moreover, the probability of (Ib) is
bounded by
\[
c_N \IP^1(\Delta\cZ[1,{n_N}]=1) { \sum_{a=1}^{n-1} \IP(\Delta\cZ[a,n-1]=1)}.
\]
The sum is the expected outdegree of vertex $n$, which, by Lemma \ref
{le0602-1}, is uniformly bounded,
and, hence, one of the events (Ia), (Ib) or (II) occurs in one step
with probability less than a constant multiple of
$c_N \IP^1(\Delta\cZ[1,{n_N}]=\nolinebreak[4]1)$.
As there are at most $c_N$ exploration steps until we end in one of the
states (A), (B) or (C),
the coupling fails due to (E1) or (E2) with a probability bounded from
above by a constant multiple of
\[
c_N^2 \IP^1(\Delta\cZ[1,{n_N}]=1) \leq 
f(1) \frac{c_N^2}{{n_N}^{1-\gamma^+}}\to0 ,
\]
in other words, the coupling succeeds with high probability.
\end{pf}

\begin{pf*}{Proof of Proposition~\ref{corcoup1}}
Apply the coupling of Lemma~\ref{premature} with $(n_N)_{N\in\IN}$ satisfying
\[
\lim_{N\to\infty} \frac{\log{n_N}}{\log N}=0 \quad\mbox{and}\quad \lim
_{N\to\infty} \frac{(\log N\log\log N)^2}{n_N^{1-\gamma^+}}=0.
\]
Then, by Lemma~\ref{initial}, we get that with high probability
%
%
\begin{equation}\label{eq0625-4}
\mbox{coupling ends in (B)} \quad\Longrightarrow\quad\#\cC_N(V)\ge c_N.
\end{equation}
As in the proof of Lemma~\ref{initial} one gets
\[
\lim_{N\to\infty} \max_{v=1,\ldots,n_N}\IP\bigl(\#\IT(v)<c_N\bigr)=0
\]
so that implication (\ref{eq0625-4}) is also valid for $\#\cC_N(V)$
replaced by $\#\IT(V)$.
Since the coupling succeeds 
we have, with high probability,
\[
\mbox{coupling ends in (A) or (B)} \quad\Longleftrightarrow\quad\#\cC
_N(V)\wedge\#\IT(V)\ge c_N,
\]
and the statement follows immediately.
\end{pf*}

\section{The idealized exploration process}\label{se7}

We now have the means to explain heuristically the approximation of the
local neighborhood of a randomly chosen vertex $V\in\mathcal G_N$
by the idealized random tree $\mathfrak T$ featuring in our main
theorems. Vertices in the network $\mathcal G_N$ are mapped onto particles
on the negative halfline in such a way that the vertex with index $n\in
\{1,\ldots,N\}$ is mapped onto position $t_n-t_N$; recall \eqref
{tndef}. Note
that the youngest vertex is placed at the origin, and older vertices
are placed to the left with decreasing intensity. In particular the
position of the particle corresponding to a vertex with fixed index
will move to the left as $N$ is increasing.

Looking at a fixed observation window $[a,b]$ on the negative halfline,
as $N\uparrow\infty$, we see that the number of particles in
the window is increasing. At the same time the age of the vertex
corresponding to a particle closest to a fixed position in the window
is increasing, which means that the probability of edges between two
such vertices is decreasing. As we shall see below, the combination of
these two effects leads to convergence of the distribution of offspring
locations on the halfline. In particular, thanks to the independence of
edges with a common right endpoint, offspring to the left converge to a
Poisson process by the law of small numbers, while offspring to the
right converge to the point processes corresponding to the pure birth
process $(Z_t \dvtx t\geq0)$ if there is no dependence on previous
generations.

The considerations of Section~\ref{se1stcoupling} suggest that the
only form of dependence of the offspring distribution
of a vertex on previous generations, is via the relative position of
its father. This information is encoded in the type of a
particle, where type $\ell$ indicates that its father is to the left of
the particle, and a numerical type $\tau$ indicates that the father is
positioned $\tau$ units to its right. It should be noted that the
\textit
{relative} positions of offspring particles only depend on the
absolute position of
the reproducing particle via the removal of particles whose position is
not in the left halfline, and which therefore do not correspond to vertices
in the network $\mathcal G_N$. This fact produces the random walk
structure, which is crucial for the analysis of the underlying tree.
Our main aim now is to prove the following result.
\begin{prop}\label{exploresurvive}
Suppose that $(c_N)_{N\iN}$ is a sequence of integers with
\[
\lim_{N\to\infty} \frac{c_N}{\log N\log\log N} =0.
\]
Then each pair $(V,\cG_N)$
can be coupled with $\mathfrak T$ such that with high probability
\[
\#\cC_N(V) \wedge c_N=\#\mathfrak{T}\wedge c_N.
\]
%
\end{prop}

We have seen so far that the neighborhood of a vertex $v$ in a large
network is similar to the random tree ${\mathbb T}(v)$ constructed
in Section~\ref{se51}. To establish the relationship between ${\mathbb T}(V)$,
for an initial vertex $V$
chosen uniformly from $\{1,\ldots,N\}$, and the idealized neighborhood
tree $\mathfrak{T}$ we apply
the projection
\[
\pi_N \dvtx(-\infty,0]\to\{1,\ldots,N\},
\]
which maps $t\leq0$ onto the smallest $m\in\{1,\ldots,N\}$ with $t\leq
-t_N+t_{m}$, to each element of the INT $\mathfrak{T}$.
We obtain a branching process with location parameters in
$\{1,\ldots,N\}$, which we call $\pi_N$-projected INT. We need to show,
using a suitable coupling, that when the INT is started with
a vertex $-X$, where $X$ is standard exponentially distributed, then
this projection is close to the random tree ${\mathbb T}(V)$.
Again we apply the concept of an exploration process.

To this end we show that, for every $v\leq0$, the $\pi_N$-projected
descendants of $v$ have a similar distribution as
the descendants of a vertex in location $\pi_N(v)$ in the labeled tree
of Section~\ref{se51}.
We provide couplings of both distributions and control the probability
of them to fail.

\subsection*{Coupling the evolution to the right for $\ell$-type vertices}
We fix $v\leq0$ and \mbox{$N\iN$}, and suppose that $m:=\pi_N(v)\geq2$. For
an $\ell$-type vertex in $v$ the cumulative sum of
$\pi_N$-projected right descendants is distributed as
$(Z_{t_n-t_N-v})_{m\leq n\leq N}$. This distribution has to be compared
with the distribution of $(\cZ[m,n])_{m\leq n\leq N}$,
which is the cumulative sum of right descendants of $m$ in ${\mathbb T}(v)$.
\begin{lemma}\label{le0604-1} 
Fix $T,N\iN$ and $v\leq0$ with $\pi_N(v)=m\in\{2,3,\ldots,N\}$. We can
couple the
processes $(Z_{t_n-t_N-v} \dvtx n\geq m)$ and $( \cZ[m,n]
\dvtx n\geq m)$ such that\vspace*{1pt}
for the coupled processes $(\cY^{\ssup1}[n] \dvtx n\geq m)$ and $(\cY
^{\ssup2}[n] \dvtx n\geq m)$ we have
\[
\IP\bigl(\cY^{\ssup1}[n]\not=\cY^{\ssup2}[n] \mbox{ for some } n\leq
\tau\bigr) \leq\bigl(f(0)+ f(T)^2\bigr) \frac1{m-1},
\]
where $\tau$ is the first time when one of the processes reaches or
exceeds $T$.
\end{lemma}
\begin{pf}
We define the process $\cY=((\cY^{\ssup1}[n],\cY^{\ssup2}[n]) \dvtx
n\geq m)$
to be the Markov process with starting distribution $\cL
(Z_{t_m-t_N-v})\otimes\delta_0$ and transition kernels $p^{\ssup n}$
such that the first and second marginal are the respective transition
probabilities of $(Z_{t_n-t_N-v} \dvtx n\geq m)$ and $(\cZ[m,n] \dvtx
n\geq m)$ and, for any integer $a\geq0$, the law $p^{\ssup n}((a,a),
\cdot)$ is the coupling of the laws
of $Z_{\Delta t_n}$ and $\cZ[n,n+1]$\vadjust{\goodbreak} under $\IP^a$ provided in
Lemma~\ref{le0525-1}. Then the processes $(\cY^{\ssup1}[n] \dvtx n\geq
m)$ and $(\cY^{\ssup2}[n] \dvtx n\geq m)$ are distributed as stated in
the lemma. Moreover, letting $\sig$ denote the first time when they disagree,
we get
\begin{eqnarray*}
\IP(\sig\leq\tau)&=& \sum_{n=m}^\infty\IP(\tau\geq n,\sig=n)\\[-3pt]
&\leq&\IP(\sig=m)+ \sum_{n=m}^\infty\IP(\sig=n+1| \tau> n,\sig
> n )
\end{eqnarray*}
and, by Lemma~\ref{le0525-1},
\[
\IP(\sig=n+1| \tau> n,\sig> n)\leq\biggl(f(T)\frac1n\biggr)^2\qquad
\mbox{for $n\in\{m,m+1,\ldots\}$.}
\]
Moreover, $\IP(\sig=m)=\IP(\cY^{\ssup1}[m]>0)=1-e^{-(t_m-t_N-v)f(0)}\leq
\frac{f(0)}{m-1}$. Consequently,
\[
\IP(\sig\leq\tau) \leq\frac{f(0)}{m-1} +f(T)^2 \sum_{n=m}^\infty
\frac{1}{n^2} \leq\bigl(f(0)+ f(T)^2\bigr) \frac1{m-1}.\vspace*{-2pt}
\]
\upqed\end{pf}

\subsection*{Coupling the evolution to the left}

Recall that a vertex $v\leq0$ produces a Poissonian number of $\pi
_N$-projected descendants at the location
$m\leq\pi_N(v)$ with parameter
%
%
\begin{equation}\label{eq0526-1}
\lam:=\int_{-t_N+t_{m-1}}^{(-t_N+t_m)\wedge v} e^{-(v-u)} \IE
[f(Z_{v-u})] \dd u.
\end{equation}
Here we adopt the convention that $t_0=-\infty$. A vertex in location
$n:=\pi_N(v)$ in ${\mathbb T}[v]$
produces a Bernoulli distributed number of descendants in $m$ with
success probability
$\IP(\Delta\cZ[m,n-1]=1)$
for $m<n$ and success probability zero for \mbox{$m=n$}. 
The following lemma provides a coupling of both distributions.\vspace*{-2pt}
\begin{lemma}\label{le0603-1} 
There exists a constant $C_{\mmbox{\ref{le0603-1}}}>0$ such that the
following holds:
Let $m,N\iN$ and $v\leq0$ with $m\leq n:=\pi_N(v)$ and define $\lam$
as in (\ref{eq0526-1}).
If $m<n$, one can couple a $\operatorname{Poiss}(\lam)$ distributed random
variable with $\Delta\cZ[m,n-1]$,
such that the coupled random variables $\Upsilon^{\ssup1}$ and
$\Upsilon^{\ssup2}$ satisfy
\[
\IP\bigl(\Upsilon^{\ssup1}\not= \Upsilon^{\ssup2}\bigr)\leq C_{\mmbox{\ref
{le0603-1}}} \frac1{m^{1+\gamma^+}} \frac1{n^{1-\gamma^+}}.
\]
%
If $m=n$, a $\operatorname{Poiss}(\lam)$ distributed random variable
$\Upsilon
^{\ssup1}$ satisfies
\[
\IP\bigl(\Upsilon^{\ssup1}\not= 0\bigr)\leq C_{\mmbox{\ref{le0603-1}}} \frac1n.\vspace*{-2pt}
\]
\end{lemma}
\begin{pf}
It suffices to prove the second statement for $m=n\geq2$.
Note that $u\mapsto e^{-u} \IE[f(Z_u)]$ is decreasing so that
\[
\lam\leq\int_{-t_N+t_{n-1}}^{v} e^{-(v-u)} \IE[f(Z_{v-u})] \dd u
\leq f(0) \sfrac{1}{n-1},\vadjust{\goodbreak}
\]
which leads directly to the second statement of the lemma.
Next, consider the case where $2\leq m<n$. Note that for $u\in
(-t_N+t_{m-1}, -t_N+t_m]$,
one has $v-u\in(t_{n-1}-t_m, t_{n}-t_{m-1})$
which, using again that $u\mapsto e^{-u} \IE[f(Z_u)]$ is decreasing,
implies that
\[
\sfrac{1}{m-1} e^{-(t_{n}-t_{m-1})} \IE[f(Z_{t_{n}-t_{m-1}})]\leq
\lam\leq\sfrac{1}{m-1} e^{-(t_{n-1}-t_{m})} \IE[f(Z_{t_{n-1}-t_{m}})].
\]
Next, note that by definition of $t_n$ we have $\log\frac nm\leq
t_n-t_m\leq\log\frac{n-1}{m-1}$ so that
%
%
\begin{eqnarray}\label{eq0529-6}
&&\biggl(1-\sfrac{1}{m-1}\biggr) \sfrac{1}{n-1} \IE[f(Z_{t_{n-1}-t_{m}})]\nonumber\\[-8pt]\\[-8pt]
&&\qquad\leq\lam
\leq\biggl(1+\sfrac{1}{m-1}\biggr) \sfrac{1}{n-1}
\IE[f(Z_{t_{n-1}-t_{m}})].\nonumber
\end{eqnarray}
On the other hand, $\Delta\cZ[m,n-1]$ is a Bernoulli random variable
with success probability
\[
p:=\sfrac{1}{n-1} \IE\bigl[f(\cZ[m,n-1])\bigr].
\]
By Lemma~\ref{lePoiBer} it suffices to control $\lam^2$ and $|\lam
-p|$. By Proposition~\ref{SGapprox} and (\ref{eq0529-6}),
%
%
\begin{equation}\label{eq0602-1}
|\lam-p|\leq C \frac1{m-1}\frac1{n-1}
\bigl(\IE[f(Z_{t_{n-1}-t_m})]+\IE\bigl[f(\cZ[m,n-1])\bigr]\bigr)
\end{equation}
and
%
%
\begin{equation}\label{eq0602-2}
\lam^2\leq4 \biggl(\sfrac{1}{n-1}\biggr)^2 \IE[f(Z_{t_{n-1}-t_m})]^2.
\end{equation}
Since $t_{n-1}-t_m\leq\log\frac{n-2}{m-1}$,
we get with Lemmas~\ref{cor0602-1} and~\ref{le0602-1} that
\[
\IE[f(Z_{t_{n-1}-t_m})]+ \IE\bigl[f(\cZ[m,n-1])\bigr] \leq C \biggl(\frac nm
\biggr)^{\gamma^+}.
\]
Recalling that $n>m\geq2$, it is now straightforward to deduce the
statement from equations (\ref{eq0602-1}) and (\ref{eq0602-2}).
It remains to consider the case where $1=m<n$. Here, we apply Lemma \ref
{cor0602-1} and $t_{n-1}\geq\log(n-1)$ to deduce that
\begin{eqnarray*}
\lam& \leq&\int_{-\infty}^{-t_N+t_1} e^{-(v-u)} \IE[f(Z_{v-u})]
\dd u \\
&\leq&C \int_{t_{n-1}}^\infty e^{-(1-\gamma^+)u} \dd u \leq
\frac{C}{1-\gamma^+} (n-1)^{\gamma^+-1},
\end{eqnarray*}
while, by Lemma~\ref{le0602-1}, $\IP(\Delta\cZ[1,n-1]=1)\leq f(0)
(n-1)^{\gamma^+-1}$,
so that a Poiss$(\lam)$ distributed random variable can be coupled with
$\Delta\cZ[1,n-1]$ so that they
disagree with probability less than a constant multiple of $n^{\gamma^+-1}$.
\end{pf}
\begin{remark}\label{rem0603-1}
Lemma~\ref{le0603-1} provides a coupling for the mechanisms with which
both trees produce \textit{left}
descendants. Since the number of descendants in individual locations
form an independent sequence of random variables, we can apply the
coupling of the lemma sequentially for each location and obtain a
coupling of the $\pi_N$-projected left descendants of a vertex $v$ and
the left descendants of $n:=\pi_N(v)$ in ${\mathbb T}[v]$. Indeed,
under the assumptions of Lemma~\ref{le0603-1}, one finds a coupling of
both processes such that
%
\[
\IP(\mbox{left descendants disagree}) \leq C_{\mmbox{\ref{le0603-1}}}
\frac1n
+C_{\mmbox{\ref{le0603-1}}} \frac1{n^{1-\gamma^+}} \sum_{m=1}^{n-1}\frac
1{m^{1+\gamma^+}}
\leq C_{\mmbox{\ref{rem0603-1}}} \frac1{n^{1-\gamma^+}},
\]
where $C_{\mmbox{\ref{rem0603-1}}}$ is a suitable positive constant.
\end{remark}

\subsection*{Coupling the evolution to the right for particles of
type $\tau\not=\ell$}

We fix $v\leq0$ and $N\iN$, and suppose that $m:=\pi_N(v)\geq2$. Also
fix a type $\tau<-v$ with $l:=\pi_N(v+\tau)>m$. The cumulative sum of
$\pi_N$-projected right descendants of a vertex $v$ of type $\tau$
(including its predecessor) is distributed according to
$(Z_{-t_N+t_n-v} \dvtx m\leq n\leq N)$
conditioned on $\Delta Z_{\tau}=1$. The cumulative sum of right
descendants in ${\mathbb T}[v]$ of a vertex in $m$ of type $l$
(including the predecessor) is distributed according to the law of $(\cZ
[m,n] \dvtx m\leq n\leq N)$ conditioned on $\Delta\cZ[m,l-1]=1$. Both
processes are Markov processes and we provide a coupling of their
transition probabilities.
\begin{lemma}\label{le0604-3} 
There exists a constant $C_{\mmbox{\ref{le0604-3}}}>0$ such that the
following holds:
Let $k\geq0$, $m,n\geq1$ be integers with $k+1<m<n$, and let %
$\tau\in(t_{n}-t_{m}, t_{n+1}-t_{m}]$. Then the random variables
$Z_{\Delta t_m}$ under $\IP^k( \ccdott  \Delta Z_{\tau}=1)$
and $\cZ[m,m+1]$ under $\IP^k( \ccdott  \Delta\cZ[m,n]=1)$ can be
coupled such that the resulting random variables
$\Upsilon^{\ssup1}$ and $\Upsilon^{\ssup2}$ satisfy
\[
\IP\bigl(\Upsilon^{\ssup1}\not=\Upsilon^{\ssup2}\bigr)\leq
C_{\mmbox{\ref{le0604-3}}}
\biggl(\frac{f(k)}{m}\biggr)^2.
\]
\end{lemma}
\begin{pf}
We couple $\Upsilon^{\ssup1}$ and $\Upsilon^{\ssup2}$ by plugging a
uniform random variable on $(0,1)$ in the generalised inverses of the
respective distribution functions and conclude that
%
\[
\IP\bigl(\Upsilon^{\ssup1}\not= \Upsilon^{\ssup2}\bigr)= \bigl|\IP\bigl(\Upsilon^{\ssup
1}=k\bigr)-\IP\bigl(\Upsilon^{\ssup2}=k\bigr)\bigr| + \IP\bigl(\Upsilon^{\ssup1}\geq k+2\bigr).
\]
The second error term is of the required order since, by Lemma~\ref{ledom},
\[
\IP\bigl(\Upsilon^{\ssup1}\geq k+2\bigr)\leq\IP^{k+1}(Z_{1/m}\geq k+3)\leq
\biggl(\frac{f(k+2)}m\biggr)^2.
\]
It remains to analyze the first error term. We have
\[
\IP\bigl(\Upsilon^{\ssup2}=k\bigr)=1-f(k) \Delta t_m \frac{P_{{m+1},n}
f(k+1)}{P_{m,n}f(k)}
\]
and, representing $(Z^{\tsupp\tau}_t \dvtx t\geq0)$ by its compensator,
\[
\IP\bigl(\Upsilon^{\ssup1} =k\bigr) =\exp\biggl\{-f(k) \int_{0}^{\Delta t_m} \frac
{ P_{\tau-u} f(k+1)}{ P_{\tau-u} f(k)} \dd u\biggr\}.
\]
We need to compare
\[
\frac{P_{{m+1},n} f(k+1)}{P_{m,n}f(k)} \quad\mbox{and}\quad \frac{ P_{u}
f(k+1)}{ P_{u} f(k)}\qquad
\mbox{for $u\in[t_{n}-t_{m+1}, t_{n+1}-t_m]$.}
\]
%
By Lemma~\ref{cor0602-1} and Proposition~\ref{SGapprox}, one has, for
$a\in\{k,k+1\}$ and sufficiently large $m$,
\begin{eqnarray*}
P_u f(a)&\leq&P_{t_{n+1}-t_m}f(a) \leq e^{\gamma^+ (1/m+1/n)}
P_{t_{n}-t_{m+1}}f(a)\\
&\leq&e^{\gamma^+ (1/m+1/n)} \biggl(1+C_{\mmbox{\ref{SGapprox}}}
\sfrac
{f(a)}{m}\biggr) P_{m+1,n}f(a)\\
& \leq&e^{\gamma^+ (1/m+1/n)+C_{\mmmbox{\ref{SGapprox}}}({f(a)}/{m})
} P_{m+1,n}f(a).
\end{eqnarray*}
Conversely,
\begin{eqnarray*}
P_u f(a)&\geq& P_{t_n-t_{m+1}} f(a)\geq e^{-\gamma^+ /m}
P_{t_n-t_{m}} f(a) \\
&\geq&e^{-\gamma^+ /m} \biggl(1- C_{\mmbox{\ref{SGapprox}}} \sfrac
{f(a)}{m}\biggr) P_{m,n}f(a).
\end{eqnarray*}
We only need to consider large $m$ and we may assume that $C_{\mmbox
{\ref
{SGapprox}}} \frac{f(k+1)}{m}\leq\frac12$, as otherwise we may choose
$C_{\mmbox{\ref{le0604-3}}}$ large
to ensure that the right-hand side in the display of the lemma exceeds
one. Then
\[
P_u f(a)\geq e^{-\gamma^+ /m- 2C_{\mmmbox{\ref{SGapprox}}}(
{f(a)}/{m})} P_{m,n}f(a)
\]
since $e^{-2y} \leq1-y$ for $y\in[0,1/2]$. Consequently,
\begin{eqnarray*}
&&e^{-\gamma^+ (2/m+1/n)-3C_{\mmmbox{\ref{SGapprox}}}({f(k+1)}/{m})}
\frac{P_{m+1,n}f(k+1)} {P_{m,n}f(k)} \\
&&\qquad\leq\frac{P_{u}
f(k+1)}{P_{u}f(k)} \leq e^{\gamma^+ (2/m+1/n)+3C_{\mmmbox{\ref
{SGapprox}}}
({f(k+1)}/{m})} \frac{P_{m+1,n}f(k+1)} {P_{m,n}f(k)}.
\end{eqnarray*}
Recall that, by Lemma~\ref{comp2}, $\frac{P_{m+1,n}f(k+1)}
{P_{m,n}f(k)}$ is uniformly bounded over all $k$ so that we arrive at
\[
\frac{P_{m+1,n}f(k+1)} {P_{m,n}f(k)} -C \frac{f(k)}m\leq\frac{P_{u}
f(k+1)}{P_{u}f(k)}\leq\frac{P_{m+1,n}f(k+1)} {P_{m,n}f(k)} +C \frac{f(k)}m
\]
for an appropriate constant $C>0$. Therefore,
\begin{eqnarray*}
&&\IP\bigl(\Upsilon^{\ssup1} =k\bigr)-\IP\bigl(\Upsilon^{\ssup2}=k\bigr)\\
&&\qquad\leq1\wedge\exp\biggl\{-f(k) \Delta t_m \biggl(\sfrac{P_{m+1,n}f(k+1)}
{P_{m,n}f(k)} -C \sfrac{f(k)}{m}\biggr)\biggr\}\\
&&\qquad\quad{} - \biggl(1-f(k) \Delta t_m \sfrac{P_{m+1,n}f(k+1)}
{P_{m,n}f(k)} \biggr)\\
&&\qquad\leq C \biggl(\sfrac{f(k)}{m}\biggr)^2+\frac12 \biggl(f(k) \Delta t_m
\biggl(\sfrac{P_{m+1,n}f(k+1)} {P_{m,n}f(k)} -C \sfrac{f(k)}{m}
\biggr)\biggr)^2\\
&&\qquad\leq C_{\mmbox{\ref{le0604-3}}} \biggl(\sfrac{f(k)}{m}\biggr)^2.
\end{eqnarray*}
Similarly, one finds that
\[
\IP\bigl(\Upsilon^{\ssup2}=k\bigr)-\IP\bigl(\Upsilon^{\ssup1}=k\bigr)\leq C_{\mmbox{\ref
{le0604-3}}} \biggl(\sfrac{f(k)}{m}\biggr)^2
\]
and putting everything together yields the assertion.
\end{pf}

From Lemma~\ref{le0604-3} we get the following analog of Lemma~\ref{le0604-1}.
\begin{lemma}\label{le0604-2} 
Fix a level $T\iN$. For any $v\leq0$ and $\tau\le-v$ with $\pi
_N(v)=m\in\{2,3,\ldots,N\}$ and $m<l:=\pi_N(v+\tau)$ we
can couple the processes\break $(Z_{t_n-t_N-v} \dvtx n\geq m)$ conditioned
on $\Delta Z_{\tau}=1$ and $(\cZ[m,n] \dvtx n\geq m)$ conditioned on
$\Delta\cZ[m,l-1]=1$ such that the coupled processes $(\cY^{\ssup1}[n]
\dvtx n\geq m)$ and $(\cY^{\ssup2}[n] \dvtx n\geq m)$ satisfy
\[
\IP\bigl(\cY^{\ssup1}[n]\not=\cY^{\ssup2}[n] \mbox{ for some } n\leq\sigma
\bigr) \leq C_{\mmbox{\ref{le0604-2}}} \bigl(f(T)^2+1\bigr) \frac1{m},
\]
where $\sigma$ is the first time when one of the processes reaches or
exceeds level~$T$.
\end{lemma}
\begin{pf}
We define the process $\cY=((\cY^{\ssup1}[n],\cY^{\ssup2}[n]) \dvtx
n\geq m)$ to be the
Markov process with starting distribution $\cL(Z_{t_m-t_N-v}| \Delta
Z_{\tau}=1)\otimes\delta_0$ and transition kernels $p^{\ssup n}$ such
that the first and second marginals are the conditioned transition
probabilities of $(Z_{t_n-t_N-v} \dvtx n\geq m)$ and $(\cZ[m,n] \dvtx
n\geq m)$ as stated in the lemma. In the case where $n<l-1$, we demand
that, for any integer $a\ge0$, the law
$p^{\ssup n}((a,a), \cdot)$ is the coupling of the laws of\vspace*{1pt}
$Z_{\Delta t_n}$ under $\IP^a( \ccdott  \Delta Z_{\tau
-(t_n-t_N-v)}=1)$ and
$\cZ[n,n+1]$ under $\IP^a( \ccdott  \Delta\cZ[n,l-1]=1)$ provided
in Lemma~\ref{le0604-3}. Conversely, we apply the unconditioned
coupling of Lemma~\ref{le0604-1} for $n\ge l$.
Letting $\varrho$ denote the first time when both evolutions disagree,
we get
\begin{eqnarray*}
\IP(\varrho\leq\sig)&=&\sum_{n=m}^\infty\IP(\sig\geq n,\varrho=n)\\
&\leq&
\IP(\varrho=m)+ \sum_{n=m}^\infty\IP(\varrho=n+1| \sigma>
n,\varrho> n )
\end{eqnarray*}
and, by Lemmas~\ref{le0604-1} and~\ref{le0604-3},
\[
\IP(\varrho=n+1| \sigma> n,\varrho> n)\leq C_{\mmbox{\ref{le0604-3}}}
\biggl(\frac{f(T)}n\biggr)^2 \!\!\qquad\mbox{for $n\in\{m,m+1,\ldots\}\setminus\{
l-1\}$.}
\]
Moreover, $\IP(\varrho=m)\leq\IP
^1(Z_{t_m-t_N-v}>0)=1-e^{-(t_m-v)f(1)}\leq\frac{f(1)}{m-1}$ and $\IP
(\varrho=l|\varrho\ge l, \sigma\ge l)\le\IP^T( Z_{\Delta
t_{l-1}}>T)\le f(T)\frac1m$. Consequently,
\[
\IP(\varrho\leq\sigma) \leq\frac{f(1)}{m-1}+\frac{f(T)}{m} + C_{\mmbox
{\ref
{le0604-3}}} f(T)^2 \sum_{n=m}^\infty\frac{1}{n^2}
\leq C_{\mmbox{\ref{le0604-2}}} \bigl(f(T)^2+1\bigr) \frac1{m}.
\]
\upqed\end{pf}

%
\begin{pf*}{Proof of Proposition~\ref{exploresurvive}}
We couple the labeled tree ${\mathbb T}(V)$
and the $\pi_N$-projected INT, starting with a coupling of the position
of the initial vertex $V$ and $\pi_N(-X)$, which
fails with probability going to zero, by Lemma~\ref{leexpunif}.

Again we apply the concept of an exploration process.
As before we categorise vertices as \textit{veiled}, if they have not yet
been discovered, \textit{active}, if they have been discovered, but if
their descendants have not yet been explored, and \textit{dead}, if they
have been discovered and all their descendants have been explored. In
one exploration step the leftmost active vertex is picked and its
descendants are explored in increasing
order with respect to the location parameter. We stop immediately once
one of the events (A), (B) or (C) happens. Note that in that case the
exploration of the last vertex might not be completed. Moreover, when
coupling two explorations, we also stop in the adverse event (E) that
the explored graphs disagree.
In event (B), the parameters $(n_N)_{N\iN}$ are chosen such that
\[
\lim_{N\to\infty} \frac{(\log N \log\log N)^\alpha}{n_N}=0
\quad\mbox{and}\quad \lim_{N\to\infty} \frac{\log n_N}{\log N}=0
\]
for $\alpha:=(1-\gamma^+)^{-1} \vee3$.
Noting that we never need to explore more than $c_N$
vertices, we see from Lemma~\ref{le0604-1}, Remark~\ref{rem0603-1} and
Lemma~\ref{le0604-2}
that the probability of a failure of this coupling is bounded by a constant
multiple of
\[
c_N \bigl(1+ f(c_N)^2\bigr) \frac1{{n_N}} + c_N \frac1{{n_N}^{1-\gamma^+}}
\leq\frac{c_N^3}{{n_N}} + \frac{c_N}{{n_N}^{1-\gamma^+}}
\longrightarrow0.
\]
%
Consequently, the coupling succeeds with high probability. As in
Lemma~\ref{initial} it is easy to see
that, with high probability, event (B) implies that
\[
\# \IT(V) \geq c_N \quad\mbox{and}\quad \#\mathfrak{T}\geq c_N.
\]
Hence we have 
\[
\# \IT(V)\wedge c_N = \# \mathfrak{T} \wedge c_N \qquad\mbox{with
high probability,}
\]
and the statement follows by combining this with Proposition~\ref{corcoup1}.
\end{pf*}
%


\section{The variance of the number of vertices in large clusters}\label{se6}

In this section we provide the second moment estimate needed to show
that our key
empirical quantity, the number of vertices in connected components of a
given size,
concentrate asymptotically near their mean.
\begin{prop}\label{prop0825-1}
Suppose that $(c_N)_{N\in\IN}$ and $(n_N)_{N\in\IN}$ are sequences of
integers satisfying $1\leq c_N,n_N\leq N$
such that $c_N^2 n_N^{\gamma^+-1}$ is bounded from above. Then, for a
constant $C_{\mmbox{\ref{prop0825-1}}}>0$
depending on these sequences and on $f$, we have
\begin{eqnarray*}
&&\var\Biggl(\frac1N\sum_{v=1}^N \ind{\{\#{\mathcal C}_N(v)\geq c_N\}
}\Biggr)\\
&&\qquad \leq2 \IP\bigl(\#{\mathcal C}_N(V)< c_N\mbox{ and }{\mathcal
C}_N(V)\cap\{1,\ldots,n_N\}\not=\varnothing\bigr)
\\
&&\qquad\quad{}+ \frac{c_N}{N} + C_{\mmbox{\ref{prop0825-1}}}
\frac{c_N^2}{n_N^{1-\gamma^+}},
\end{eqnarray*}
where $V$ is independent of ${\mathcal G}_N$ and uniformly distributed
on $\{1,\ldots,N\}$.
\end{prop}
\begin{pf}
Let $v, w$ be two distinct vertices of $\cG_N$. We start by exploring
the neighborhood of $v$ similarly as in Section~\ref{se1stcoupling}.
As before we classify the vertices as veiled, active and dead, and in
the beginning only $v$ is active and the remaining vertices are veiled.
In one exploration step we pick the leftmost active vertex and
consecutively (from the left to the right) explore its immediate
neighbors in the set of veiled vertices only. Newly found vertices are
activated and the vertex to be explored is set to dead after the
exploration. We immediately stop the exploration once one of the events:
\begin{longlist}[(C)]
\item[(A)] the number of unveiled vertices in the cluster reaches $c_N$,
\item[(B)] one vertex in $\{1,\ldots,n_N\}$ is activated or
\item[(C)] there are no more active vertices left,
\end{longlist}
happens. Note that when we stop due to (A) or (B) the exploration of
the last vertex might not be finished. In that case we call this vertex
\textit{semi-active}.

We proceed with a second exploration process, namely the exploration of
the cluster of $w$. This exploration follows the same rules
as the first exploration process, treating vertices that remained
active or semi-active at the end of the first exploration as veiled.
In addition to the stopping in the cases (A), (B), (C) we also stop the
exploration once a vertex is unveiled which
was also unveiled in the first exploration, calling this event (D).
We consider the following events:
\begin{longlist}[$E_{3}^{v,w}$:]
\item[$E^v$:] the first exploration started with
vertex $v$ ends in (A) or (B);\vspace*{1pt}
\item[$E_{1}^{v,w}$:] $w$ is unveiled during the first exploration
(that of $v$);\vspace*{2pt}
\item[$E_{2}^{v,w}$:] $w$ remains veiled in the first exploration and
the second exploration ends in (A) or (B)
but not in (D);
\item[$E_{3}^{v,w}$:] $w$ remains veiled in the first exploration and
the second exploration ends in (D).
\end{longlist}
%

We have
%
%
\begin{eqnarray}\label{eq0610-1}
&&\sum_{v=1}^N \sum_{w=1}^N \IP\bigl(\# {\mathcal C}_N(v)\geq c_N , \#
{\mathcal C}_N(w)\geq c_N \bigr)\nonumber\\
&&\qquad\leq
\sum_{v=1}^N \sum_{w=1}^N \sum_{k=1}^3 \IP(E^v\cap E_{k}^{v,w}) \\
&&\qquad= \sum_{v=1}^N \IP(E^v) \sum_{k=1}^3 \sum_{w=1}^N \IP(E_{k}^{v,w}\mid
E^v).\nonumber
\end{eqnarray}
As the first exploration immediately stops once one has unveiled $c_N$
vertices, we conclude that,
for fixed $v$,
%
%
\begin{equation}\label{eq0610-2}
\sum_{w=1}^N \IP(E_{1}^{v,w}\mid E^v)=\IE\Biggl[\sum_{w=1}^N
\ind_{E_{1}^{v,w}}\Big| E^v\Biggr]\leq c_N .
\end{equation}
To analyze the remaining terms, we fix distinct vertices $v$ and $w$
and note that the configuration after the first
exploration can be formally described by an element $\mathfrak k$ of
\[
\{\mbox{open, closed, unexplored}\}^{E_N},
\]
where $E_N:=\{(a,b)\in\{1,\ldots,N\}^2\dvtx i<j\}$ denotes the set of
possible edges. We pick a feasible configuration $\mathfrak k$ and
denote by $\cE_{\mathfrak k}$ the event that the first exploration
ended in this configuration. On the event $\cE_{\mathfrak k}$ the
status of each vertex (veiled, active, semi-active or dead) at the end
of the first exploration is determined.
Suppose $\mathfrak k$ is such that $w$ remained veiled in the first
exploration, which means that $\cE_{\mathfrak k}$ and $E_{1}^{v,w}$ are
disjoint events. Next, we note that
%
%
\begin{equation}\label{eq0610-3}
\IP(E^{v,w}_{2}| \cE_{\mathfrak k}) \leq
\IP(E^w).
\end{equation}
Indeed, if in the exploration of $w$ we encounter an edge which is open
in the configuration $\mathfrak k$, we have
unveiled a vertex which was also unveiled in the exploration of $v$,
the second exploration ends in (D) and
hence $E^{v,w}_{2}$ does not happen.
Otherwise, the event $\cE_{\mathfrak k}$ influences the exploration
of $w$ only in the sense that in the
degree evolution of some vertices some edges may be conditioned to be closed.
By Lemma~\ref{extra} this conditional probability
is bounded by the unconditional probability, and hence we obtain (\ref
{eq0610-3}).

Finally, we analyze the probability $\IP(E^{v,w}_{3}| \cE_{\mathfrak
k})$. If the second exploration process ends in state (D),
we have discovered an edge connecting the exploration started in $w$ to\vadjust{\goodbreak}
an active or semi-active vertex $a$ from the first exploration.
Recall that in each exploration we explore the immediate neighborhoods
of at most $c_N$ vertices. Let $\mathfrak K\in E_N$ be a
feasible configuration at the beginning of the neighborhood exploration
of a vertex $n>n_N$ and note that this implies every edge
which is open (resp., closed) in $\mathfrak k$ is also open (resp.,
closed) in $\mathfrak K$. Recall that $\cE_{\mathfrak K}$ denotes the
event that this configuration
is seen in the combined exploration processes. We denote by $\mathfrak
a$ and $\mathfrak s$ the set of active and semi-active vertices of the
\textit{first exploration} induced by $\mathfrak k$ (or, equivalently, by
$\mathfrak K$). Moreover, we denote by $\mathfrak d$ the set of
dead vertices of the combined exploration excluding the father of $n$,
and, for $a\in\mathfrak a\cup\mathfrak s$, we let $\mathfrak d_a$
denote the set of dead vertices of the ongoing exploration excluding
the father of $n$, plus the vertices that were marked as dead in the
first exploration at the time the vertex $a$ was discovered.
We need to distinguish several cases.

\textit{First}, consider the case $a\in\mathfrak a$ with $a<n$. By
definition of the combined exploration process, we know that $a$ has no
jumps in its indegree evolution at times associated to the vertices
$\mathfrak d_a$.
If $a$ was explored from the right, say with father in $b$, we thus get
%
%
\begin{eqnarray}\label{eq0916-1}
&&\IP(\exists\mbox{ edge between $a$ and $n$}| \cE_{\mathfrak
K}) \nonumber\\
&&\qquad= \IP(\Delta\cZ[a,n-1]=1
|\Delta\cZ[a,b-1]=1\mbox{ and } \Delta\cZ[a, d-1]=0\\
&&\hspace*{273pt} \forall d\in
\mathfrak d_a).
\nonumber
\end{eqnarray}
If $a$ was explored from the left, then
%
%
\begin{eqnarray}\label{eq0916-2}
&&\IP(\exists\mbox{ edge between $a$ and $n$}| \cE_{\mathfrak
K})\nonumber\\[-8pt]\\[-8pt]
&&\qquad= \IP(\Delta\cZ[a,n-1]=1| \Delta\cZ[a, d-1]=0\
\forall d\in\mathfrak
d_a).\nonumber
\end{eqnarray}
\textit{Second}, consider the case $a\in\mathfrak a$ with $n<a$.
By definition of the combined exploration process, the indegree
evolution of $n$ has no jumps that can be associated to edges
connecting to $\mathfrak d$.

Hence, if $n$ was explored from the right, say with father in $b$, then
%
%
\begin{eqnarray}
\label{eq0916-3}\quad
&&\IP(\exists\mbox{ edge between $a$ and $n$}| \cE_{\mathfrak
K}) \nonumber\\
&&\qquad= \IP(\Delta\cZ[n,a-1]=1
| \Delta\cZ[n,b-1]=1\mbox{ and } \Delta\cZ[n, d-1]=0\\
&&\hspace*{278pt}\forall
d\in\mathfrak d),
\nonumber
\end{eqnarray}
and, if $n$ was explored from the left, then
%
%
\begin{eqnarray}\label{eq0916-4}
&&\IP(\exists\mbox{ edge between $a$ and $n$}| \cE_{\mathfrak
K})\nonumber\\[-8pt]\\[-8pt]
&&\qquad= \IP(\Delta\cZ[n,a-1]=1| \Delta\cZ[n, d-1]=0\
\forall d\in\mathfrak
d).\nonumber
\end{eqnarray}
\textit{Third}, consider $a\in\mathfrak s$ and denote by $a'$ the last
vertex which was unveiled in the first exploration. If $a'>n$, then the
existence of an edge between $a$ and $n$ was already explored in the
first exploration,\vadjust{\goodbreak} and no edge was found. If $a'<n<a$,
we find estimates \eqref{eq0916-3}, \eqref{eq0916-4} again. If $a<n$
and the father $b$ of $a$ satisfies $b>a'\vee a$,
%
%
\begin{eqnarray}
\label{eq0916-5}
&&\IP(\exists\mbox{ edge between $a$ and $n$} | \cE_{\mathfrak
K}) \nonumber\\
&&\qquad\leq
\sup_{0\le k \le c_N-1} \IP^k(\Delta
\cZ[a\vee a',n-1]=1 | \Delta\cZ[a\vee a',b-1]=1\\
&&\hspace*{135.4pt}\mbox{ and } \Delta
\cZ[a \vee a', d-1]=0\
\forall d\in\mathfrak d_a),
\nonumber
\end{eqnarray}
and if $a=v$ or the father $b$ of $a\vee a'$ satisfies $b< a \vee a'$,
%
%
\begin{eqnarray}
\label{eq0916-6}\qquad
&&\IP(\exists\mbox{ edge between $a$ and $n$}| \cE_{\mathfrak
K})\nonumber\\[-8pt]\\[-8pt]
&&\qquad\leq
\sup_{0\le k \le c_N}\IP^k(\Delta\cZ[a \vee a',n-1]=1| \Delta
\cZ[a \vee a', d-1]=0\ \forall d\in\mathfrak d_a).
\nonumber
\end{eqnarray}

Using first Lemma~\ref{condaway}, then Lemmas~\ref{ledom2} and
\ref{le0416-1} we see that
the terms \eqref{eq0916-1}--\eqref{eq0916-4} are bounded by
\[
C_{\mmbox{\ref{condaway}}} \IP^1(\Delta\cZ[a,n-1]=1) \leq C_{\mmbox{\ref
{condaway}}} \frac{P_{1,n_N}f(1)}{n_N},
\]
and similarly, the terms \eqref{eq0916-5}--\eqref{eq0916-6} are bounded by
\[
C_{\mmbox{\ref{condaway}}} \IP^{c_N}(\Delta\cZ[a,n-1]=1) \leq
C_{\mmbox{\ref{condaway}}} \frac{P_{1,n_N}f(c_N)}{n_N}.
\]


Note that there are at most $c_N$ vertices $a\in\mathfrak a\cup
\mathfrak s$ and at most one
of those is semi-active. For each of these $a$ we have to test the
existence of edges no more than
$c_N$ times. Hence, using also Lemma~\ref{le0602-1} and the boundedness
of $f(n)/n$, we find
$C_{\mmbox{\ref{prop0825-1}}}>0$ such that
\begin{eqnarray*}
\IP(E_{3}^{v,w} | E^v ) &\leq& C_{\mmbox{\ref{condaway}}} c_N^2
\frac{P_{1,n_N}f(1)}{n_N}
+ C_{\mmbox{\ref{condaway}}} c_N \frac{P_{1,n_N}f(c_N)}{n_N}\\
&\leq& C_{\mmbox{\ref{prop0825-1}}} \frac{c_N^2}{n_N^{1-\gamma^+}}.
\end{eqnarray*}
Summarizing our steps, we have
\begin{eqnarray*}
&& \var\Biggl( \frac1N \sum_{v=1}^N \ind{\{\#{\mathcal C}_N(v) \geq c_N\}
}\Biggr) \\
&&\qquad \leq\IE\Biggl[ \frac1{N^2}\sum_{v=1}^N \sum_{w=1}^N \ind\{ \#\cC
_N(v) \geq c_N, \#\cC_N(w) \geq c_N\}\Biggr]\\
&&\qquad\quad{} - \frac1{N^2} \sum_{v=1}^N \sum_{w=1}^N \IP(E^v)\IP
(E^w) \\[-2pt]
&&\qquad\quad{} + 2 \frac1{N} \sum_{v=1}^N \IP\bigl(\#{\mathcal C}_N(v)<
c_N\mbox{ and }{\mathcal C}_N(v)\cap\{1,\ldots,n_N\}\not=\varnothing\bigr)
\\[-2pt]
&&\qquad\leq2 \IP\bigl(\#{\mathcal C}_N(V)< c_N\mbox{ and }{\mathcal
C}_N(V)\cap\{1,\ldots,n_N\}\not=\varnothing\bigr)\\
&&\qquad\quad{}
+ \frac{c_N}{N} + C_{\mmbox{\ref{prop0825-1}}} \frac
{c_N^2}{n_N^{1-\gamma^+}}
\end{eqnarray*}
as required to complete the proof.
\end{pf}

\section{\texorpdfstring{Proof of Theorem \protect\ref{main2}}{Proof of Theorem 1.8}}\label{se8}
We start by proving the lower bound for $\cC^{\ssup1}_N$. Suppose
therefore that $p(f)>0$,
fix $\delta>0$ arbitrarily small and use Lem\-ma~\ref{pcontinuous} to
choose $\eps>0$ such that the survival
probability of $\bar f=f-\eps$ is larger than $p(f)-\delta$.
We denote by $(\bar\cG_N)_{N\iN}$ a sequence of random networks with
attachment rule $\bar f$ and let
$\bar{\mathcal C}_N(v)$ the connected component of $v$ in $\bar\cG
_N$. Suppose a vertex $V$ is
chosen uniformly at random from $\{1,\ldots,N\}$.
We choose $c_N := \lfloor\log N \sqrt{\log\log N} \rfloor$ and
observe that, by Proposition~\ref{exploresurvive},
%
%
\begin{eqnarray}\label{eq0630-4}
&&\IE\Biggl[\frac1{N} \sum_{v=1}^N \ind{\{\#\bar{\mathcal C}_N(v)\geq c_N\}
}\Biggr]\nonumber\\[-9pt]\\[-9pt]
&&\qquad= \IP\{\#\bar{\mathcal C}_N(V)\geq c_N\}
\quad\longrightarrow\quad \IP\{\#\mathfrak T=\infty\}
\geq
p(f)-\delta\nonumber
\end{eqnarray}
as $N$ tends to infinity.
By Proposition~\ref{prop0825-1}
with $n_N:=\lfloor(\log N)^{4/({1-\gamma^+})} \rfloor$, we have
\begin{eqnarray*}
&&\var\Biggl(\frac1{N}\sum_{v=1}^N \ind{\{\#\bar{\mathcal C}_N(v)\geq
c_N\}}\Biggr)\\[-2pt]
&&\qquad\leq2 \IP\bigl(\#\bar{\mathcal C}_N(V)< c_N \mbox{ and }\bar
{\mathcal C}_N(V)\cap\{1,\ldots,{n_N}\}\not=\varnothing\bigr)\\[-2pt]
&&\qquad\quad{}+ \frac{c_N}{N} + C_{\mmbox{\ref{prop0825-1}}} \frac
{c_N^2}{n_N^{1-\gamma^+}}.
\end{eqnarray*}
The first summand goes to zero by Lemma~\ref{initial}, and so do the
remaining terms by the choice of
our parameters. Hence
\[
\liminf_{N \to\infty} \frac1{N} \sum_{v=1}^N \ind{\{\#\bar{\mathcal
C}_N(v)\geq c_N\}} \geq p(f)-\delta\qquad
\mbox{in probability,}
\]
and Proposition~\ref{sprinkling} implies that, with high probability,
there exists a connected component comprising
at least a proportion $p(f)$ of all vertices, proving the lower bound.

To see the upper bound we work with the original attachment
function~$f$. In analogy to (\ref{eq0630-4}) we obtain
\[
\lim_{N\to\infty} \IE\Biggl[\frac1{N} \sum_{v=1}^N \ind{\{\#{\mathcal
C}_N(v)\geq c_N\}}\Biggr] =p(f).\vadjust{\goodbreak}
%
\]
%
As in the lower bound, the variance goes to zero, and hence we have
\[
\lim_{N \to\infty} \frac1{N} \sum_{v=1}^N \ind{\{\#{\mathcal
C}_N(v)\geq c_N\}} = p(f)\qquad
\mbox{in probability.}
\]
From this we infer that, in probability,
\[
\limsup_{N \to\infty} \frac{\# {\mathcal C}_N^{\ssup1}}{N} \leq
\limsup_{N \to\infty} \frac{c_N}N \vee\Biggl( \frac1{N} \sum_{v=1}^N
\ind{\{\#{\mathcal C}_N(v)\geq c_N\}} \Biggr)
\leq p(f)
\]
proving the upper bound.

Finally, to prove the result on the size of the second largest
connected component, note that we have seen in particular
that
\[
\lim_{N \to\infty} \frac1{N} \sum_{v=1}^N \ind{\{\#{\mathcal
C}_N(v)\geq c_N\}} = p(f)\qquad
\mbox{in probability,}
\]
so that, with high probability, the proportion of vertices in clusters
of size greater or equal $c_N$ is asymptotically
equal to the proportion of vertices in the giant component. This
implies that the proportion of vertices, which are
not in the giant component but in components of size at least $c_N$
goes to zero in probability, which is a
stronger result than the stated claim.

\section{\texorpdfstring{Proof of Theorem \protect\ref{main3}}{Proof of Theorem 1.9}}\label{se9}

We fix $k\iN$ and choose $c_N:=k+1$.
By Proposition~\ref{exploresurvive}, we have
\[
\lim_{N \to\infty} \frac1{N} \IE\Biggl[\sum_{v=1}^N \ind{\{\#
{\mathcal
C}_N(v)\leq k\}}\Biggr] =\lim_{N\to\infty} \IP\bigl(\#{\mathcal C}_N(V)\leq
k\bigr)= \IP(\#\mathfrak T\leq k),
\]
and Proposition~\ref{prop0825-1} yields
\[
\var\biggl(\frac1N \ind{\{\#{\mathcal C}_N(v)\leq k\}}\biggr)
=\var\biggl(\frac1N \ind{\{\#{\mathcal C}_N(v) \geq c_N\}}\biggr)\to0.
\]
This implies the statement, as $k$ is arbitrary.

\section{\texorpdfstring{Proof of Theorem \protect\ref{percol}}{Proof of Theorem 1.6}}\label{se10}

The equivalence of the divergence of the sequence in Theorem \ref
{percol}, and the criterion
$\mathcal I=\varnothing$ stated in (i) of Remark~\ref{percol2}
follows from the bounds on the spectral radius of
the operators $A_\alpha$ given in the proof of Proposition~\ref{main4}.
Moreover, it is easy to see from the arguments of Section~\ref{se3}
that the survival of the INT under
percolation with retention parameter $p$ is equivalent to
the existence of $0<\alpha<1$ such that
\[
\rho(pA_\alpha)=p\rho(A_\alpha)\leq1.
\]
Hence, to complete the proof of Theorem~\ref{percol} and Remark \ref
{percol2}, it suffices to
show that, for a fixed retention parameter $0<p<1$, the existence of a
giant component for the percolated\vadjust{\goodbreak}
network is equivalent to the survival of the INT under
percolation with retention parameter $p$. We now give a sketch of this
by showing how the corresponding arguments in
the proof of Theorem~\ref{main2} have to be modified.

As in the proof of Theorem~\ref{main2} the main part of the argument
consists of couplings of the
exploration process of the neighborhood of a vertex in the network to
increasingly simple objects.
To begin with we have to couple the exploration of vertices in the
percolated network
and the percolated labeled tree, using arguments as in Section \ref
{se1stcoupling}. We only modify the
exploration processes a little: whenever we find a new vertex, instead
of automatically declaring it active,
we declare it \textit{active} with probability $p$ and \textit{passive}
otherwise. We do this independently for each
newly found vertex. We still explore at every step the leftmost active
vertex, but we change the stopping
criterion (E1): we now stop the process when we rediscover an active
\textit{or passive} vertex. We also stop the process
when we have discovered more than $2 \frac{1-p}{p} c_N$ passive
vertices,\vspace*{1pt} calling this event (E3). All other stopping
criteria are retained literally.

By a simple application of the strong law of large numbers we see that
the probability of stopping in the event (E3)
converges to zero. The proof of Lem\-ma~\ref{premature} carries over to
our case, as it only uses that the number of
dead, active and passive vertices is bounded by a constant multiple
of $c_N$. Hence the coupling of explorations is
successful with high probability.

Similarly, the coupling of the exploration processes for the random
labeled tree and the idealized neighborhood tree
constructed in Section~\ref{se7} can be performed so that under the
assumption on the parameters
given in Proposition~\ref{exploresurvive}, we have 
\[
\#\cC_N^*(V)\wedge c_N=\#\mathfrak{T}^*\wedge c_N \qquad\mbox{with
high probability},
\]
where $\cC_N^*(v)$ denotes the connected component in the percolated
network, which contains the vertex $v$, and
$\mathfrak{T}^*$ is the percolated INT.

In order to analyze the variance of the number of vertices in large
clusters of the percolated network,
we modify the exploration processes described in the proof of
Proposition~\ref{prop0825-1}
a little: in the first exploration we activate newly unveiled vertices
with probability $p$ and declare them
passive otherwise. We always explore the neighborhood of the leftmost
active vertex and investigate its links to the
set of veiled \textit{or passive} vertices from left to right, possibly
activating a passive vertex when it is revisited.
We stop the exploration in the events (A), (B) and (C) as before and,
additionally, if the number of passive vertices
exceeds $2 \frac{1-p}{p} c_N$, calling this event (A\ppprime). As
before, the probability of stopping in (A\ppprime) goes to zero by
the strong law of large numbers.

The exploration of the second cluster follows the same rules as that of
the first, treating vertices that were left active, semi-active
or passive in the first exploration as veiled. In addition to the
stopping events (A), (A\ppprime), (B) and (C) we also stop in the event (D)
when a vertex is unveiled which was also unveiled in the first
exploration. This vertex may have been active, semi-active or passive
at the
end of the first exploration. We then introduce the event $E^v$ that
the first exploration ends in events (A), (A\ppprime) or (B), events
$E^{v,w}_1$ and $E^{v,w}_3$ as before, and event $E^{v,w}_2$ that $w$
remained veiled in the first exploration and the second exploration
ends in (A), (A\ppprime) or (B). We can write
\[
\sum_{v=1}^N \sum_{w=1}^N \IP\bigl(\# {\mathcal C}'_N(v)\geq c_N , \#
{\mathcal C}'_N(w)\geq c_N \bigr)
\leq\sum_{v=1}^N \IP(E^v) \sum_{k=1}^3 \sum_{w=1}^N \IP
(E_{k}^{v,w}\mid E^v),
\]
where ${\mathcal C}'_N(v)$ denotes the connected component of $v$ in
the percolated network. The summand corresponding to $k=1$ can be
estimated as before. For the other summands we describe the
configuration after the first exploration as an element $\mathfrak k$ of
\[
\{\mbox{open, closed, removed, unexplored}\}^{E_N},
\]
where edges corresponding to the creation of passive vertices are
considered as ``removed.'' We again obtain that
$\IP(E_2^{v,w}| {\mathcal E}_{\mathfrak k}) \leq\IP(E^w)$ using the
fact that if in the second exploration we ever encounter
an edge which is open or removed in the configuration ${\mathfrak k}$,
the second exploration ends in (D), and $E_2^{v,w}$
does not occur. Finally, the estimate of $\IP(E_3^{v,w}| {\mathcal
E}_{\mathfrak k})$ carries over to our situation as it relies
only on the fact that the number of unveiled vertices in the first
exploration is bounded by a constant multiple of $c_N$. We thus
obtain a result analogous to Proposition~\ref{prop0825-1}.

Using straightforward analogues of the results in Section~\ref{se4} we
can now show that the existence of a giant component for the percolated
network is equivalent to the survival of the INT under percolation with
retention parameter $p$ using the argument of
Section~\ref{se8}. This completes the proof of Theorem~\ref{percol}.

\begin{appendix}\label{app}

\section*{Appendix}

In this Appendix we provide two auxiliary coupling lemmas.
\begin{lemma}\label{lePoiBer} Let $\lam\geq0$ and $p\in[0,1]$,
$X^{\ssup1}$ Poisson distributed with parameter $\lam$, and $X^{\ssup
2}$ Bernoulli
distributed with parameter $p$. Then there exists a coupling of these
two random
variables such that
\[
\IP\bigl(X^{\ssup1}\not=X^{\ssup2}\bigr)\leq\lam^2+|\lam-p|.
\]
\end{lemma}
\begin{pf} We only need to consider the case where $\lam\in[0,1]$. Then
$X^{\ssup1}$ can
be coupled to a Bernoulli\vspace*{1pt} distributed random variable $X$ with
parameter $\lambda$, such that
$\IP(X^{\ssup1}\not=X)= \lam-\lam e^{-\lam} \leq\lam^2$. Moreover,
$X$ and $X^{\ssup2}$ can be coupled
such that $\IP(X\not=X^{\ssup2})=|p-\lam|$. The two facts together
imply the statement.
\end{pf}
\begin{lemma}\label{leexpunif} Let $Y$ be standard exponentially
distributed and $X$ uniformly
distributed on $\{1,\ldots,N\}$. Then $X$ and $Y$ can be coupled in
such a way that
\[
\IP\bigl(X\not=\pi_N(-Y)\bigr)\leq C_{\mmbox{\ref{leexpunif}}} \frac{\log N}N
\]
for the function $\pi_N$ defined at
the beginning of Section~\ref{se7}.
\end{lemma}
\begin{pf}
For $2 \leq k\leq N$, we have
\begin{eqnarray*}
\IP\bigl(\pi_N(-Y)=k\bigr) & = & \IP\Biggl( \sum_{j=k}^{N-1} \frac1j \leq Y < \sum
_{j=k-1}^{N-1} \frac1j\Biggr)\\
&=& \exp\Biggl\{ -\sum_{j=k}^{N-1} \frac1j \Biggr\} - \exp\Biggl\{ -\sum
_{j=k-1}^{N-1} \frac1j \Biggr\}.
\end{eqnarray*}
Since $\sum_{j=k-1}^{N-1} \frac1j\geq\log\frac{N}{k-1}$ and $e^x\leq
1+x+x^2$ for $x\in[1,2]$, we get
\begin{eqnarray*}
\IP\bigl(\pi_N(-Y)=k\bigr) &=& \exp\Biggl\{ -\sum_{j=k-1}^{N-1} \frac1j \Biggr\}
\bigl(e^{1/({k-1})}-1\bigr)\\
&\leq&\frac{k-1}{N} \biggl(\frac1{k-1}+\frac1{(k-1)^2}\biggr) \leq\frac
1N +\frac2{Nk}.
\end{eqnarray*}
Similarly, one obtains that $\IP(\pi_N(-Y)=k)\geq\frac1N -\frac2{Nk}$.
Hence we can couple the random variables so that, for a suitable
constant $C_{\mmbox{\ref{leexpunif}}}>0$,
\[
\IP\bigl(X\not=\pi_N(-Y)\bigr)\leq
\sum_{k=2}^N \biggl| \IP\bigl(\pi_N(-Y)=k\bigr)-\sfrac{1}{N} \biggr|
\leq C_{\mmbox{\ref{leexpunif}}} \frac{\log N}N.
\]
\upqed\end{pf}
\end{appendix}

\section*{Acknowledgments}

We would like to thank Maren Eckhoff, Christian M\"{o}nch, Tom Mountford
and two referees for heplful comments on the first version of this
paper.


%

\printaddresses

\end{document}